\documentclass[10pt,twoside]{amsart}

\usepackage{amssymb}
\usepackage{amscd}
\usepackage{amsmath}
\usepackage[all]{xy}
\CompileMatrices
\usepackage[backref]{hyperref}

\newtheorem{thm}{Theorem}[section]
\newtheorem{prop}[thm]{Proposition} 
\newtheorem{lem}[thm]{Lemma} 
\newtheorem{cor}[thm]{Corollary} 
\newtheorem{conj}[thm]{Conjecture} 
\newtheorem{defn}[thm]{Definition} 
\newtheorem{example}[thm]{Example} 

\newtheorem{remark}[thm]{Remark}

\newcommand{\z}{\mathbb{Z}}

\newcommand{\qnr}{ \widehat{\mathbb{Q}_p^{nr}}}

\newcommand{\kl}{[\![}
\newcommand{\kr}{]\!]}

\renewcommand{\projlim}[1] {{\lim\limits_{\stackrel{\displaystyle
\longleftarrow}{#1}}}}

\newcommand{\T}{\mathbb{ T}}
\newcommand{\cA}{\ifmmode {A^\cdot}\else${A^\cdot}$\ \fi}
\newcommand{\cB}{\ifmmode { B^\cdot}\else${B^\cdot}$\ \fi}
\newcommand{\cC}{\ifmmode { C^\cdot}\else${C^\cdot}$\ \fi}
\newcommand{\cP}{\ifmmode {P^\cdot}\else${P^\cdot}$\ \fi}
\renewcommand{\L}{\ifmmode {\mathcal{L}}\else$\mathcal{L}$\ \fi}
\newcommand{\ca}{\ifmmode { A^\cdot_0}\else${A^\cdot_0}$\ \fi}
\newcommand{\be}{\begin{equation}}
\newcommand{\ee}{\end{equation}}
\renewcommand{\d}{\mathbf{d}}
\renewcommand{\u}{\mathbf{1}}
\newcommand{\id}{\mathrm{id}}
\newcommand{\bbC}{\ifmmode {\mathbb{C}}\else$\mathbb{C}$\ \fi}
\newcommand{\bbR}{\ifmmode {\mathbb{R}}\else$\mathbb{R}$\ \fi}
\newcommand{\R}{\mbox{$\Bbb R$}}
\renewcommand{\r}{\mathrm{R}\Gamma}
\newcommand{\cone}{\mathrm{cone}}
\newcommand{\bq}{\mathbb Q}
\newcommand{\Q}{\mathbb Q}

\newcommand{\bz}{\mathbb Z}
\newcommand{\br}{\mathbb R}
\newcommand{\bc}{\mathbb C}

\renewcommand{\H}{\mathrm{H}}

\newcommand{\B}{\mathfrak B}

\renewcommand{\u}{\mathbf{1}}

\newcommand{\zp}{{\bz_p}}
\newcommand{\qp}{{\bq_p}}

\newcommand{\p}{\mathfrak p}

\renewcommand{\O}{\mathcal{ O}}

\newcommand{\triv}{\rm triv}

\newcommand{\ql}{{\bq_\ell}}

\newcommand{\C}{\mathcal{C}}
\newcommand{\La}{\ifmmode\Lambda\else$\Lambda$\fi}

\DeclareMathOperator{\Spec}{Spec}

\DeclareMathOperator{\Ind}{Ind}

\DeclareMathOperator{\ind}{ind}

\DeclareMathOperator{\Gal}{Gal}

\DeclareMathOperator{\Hom}{Hom}

\DeclareMathOperator{\cok}{cok} 

\renewcommand{\det}{\text{det}}

\def\YEAR{\year}\newcount\VOL\VOL=\YEAR\advance\VOL by-1995
\def\firstpage{1}\def\lastpage{1000}
\def\received{}\def\revised{}
\def\communicated{}

\makeatletter
\def\magnification{\afterassignment\m@g\count@}
\def\m@g{\mag=\count@\hsize6.5truein\vsize8.9truein\dimen\footins8truein}
\makeatother

\font\eightrm=cmr8
\font\caps=cmcsc10
\font\Caps=cmcsc10 scaled \magstep1   % Title                 % Theorem, Lemma
\makeatletter
\@specialpagefalse\headheight=8.5pt
\def\DocMath{}
\renewcommand{\@evenhead}{%
    \ifnum\thepage>\lastpage\rlap{\thepage}\hfill%
    \else\rlap{\thepage}\slshape\leftmark\hfill{\caps\SAuthor}\hfill\fi}%
\renewcommand{\@oddhead}{%
    \ifnum\thepage=\firstpage{\DocMath\hfill\llap{\thepage}}%
    \else{\slshape\rightmark}\hfill{\caps\STitle}\hfill\llap{\thepage}\fi}%
\makeatother

\def\TSkip{\bigskip}
\newbox\TheTitle{\obeylines\gdef\GetTitle #1
\ShortTitle  #2
\SubTitle    #3
\Author      #4
\ShortAuthor #5
\EndTitle
{\setbox\TheTitle=\vbox{\baselineskip=20pt\let\par=\cr\obeylines%
\halign{\centerline{\Caps##}\cr\noalign{\medskip}\cr#1\cr}}%
        \copy\TheTitle\TSkip\TSkip%
\def\next{#2}\ifx\next\empty\gdef\STitle{#1}\else\gdef\STitle{#2}\fi%
\def\next{#3}\ifx\next\empty%
    \else\setbox\TheTitle=\vbox{\baselineskip=20pt\let\par=\cr\obeylines%
    \halign{\centerline{\caps##} #3\cr}}\copy\TheTitle\TSkip\TSkip\fi%
\centerline{\caps #4}\TSkip\TSkip%
\def\next{#5}\ifx\next\empty\gdef\SAuthor{#4}\else\gdef\SAuthor{#5}\fi%
\ifx\received\empty\relax
    \else\centerline{\eightrm Received: \received}\fi%
\ifx\revised\empty\TSkip%
    \else\centerline{\eightrm Revised: \revised}\TSkip\fi%
\ifx\communicated\empty\relax
    \else\centerline{\eightrm Communicated by \communicated}\fi\TSkip\TSkip%
\catcode'015=5}}\def\Title{\obeylines\GetTitle}
\def\Abstract{\begingroup\narrower
    \parskip=\medskipamount\parindent=0pt{\caps Abstract. }}
\def\EndAbstract{\par\endgroup\TSkip}

\long\def\MSC#1\EndMSC{\def\arg{#1}\ifx\arg\empty\relax\else
     {\par\narrower\noindent%
     1991 Mathematics Subject Classification: #1\par}\fi}

\long\def\KEY#1\EndKEY{\def\arg{#1}\ifx\arg\empty\relax\else
        {\par\narrower\noindent Keywords and Phrases: #1\par}\fi\TSkip}

\newbox\TheAdd\def\Addresses{\vfill\copy\TheAdd\vfill
    \ifodd\number\lastpage\vfill\eject\phantom{.}\vfill\eject\fi}
{\obeylines\gdef\GetAddress #1
\Address #2
\Address #3
\Address #4
\EndAddress
{\def\xs{4.3truecm}\parindent=0pt
\setbox0=\vtop{{\obeylines\hsize=\xs#1\par}}\def\next{#2}
\ifx\next\empty % 1 address
     \setbox\TheAdd=\hbox to\hsize{\hfill\copy0\hfill}
\else\setbox1=\vtop{{\obeylines\hsize=\xs#2\par}}\def\next{#3}
\ifx\next\empty % 2 addresses
     \setbox\TheAdd=\hbox to\hsize{\hfill\copy0\hfill\copy1\hfill}
\else\setbox2=\vtop{{\obeylines\hsize=\xs#3\par}}\def\next{#4}
\ifx\next\empty\ % 3 addresses
     \setbox\TheAdd=\vtop{\hbox to\hsize{\hfill\copy0\hfill\copy1\hfill}
                \vskip20pt\hbox to\hsize{\hfill\copy2\hfill}}
\else\setbox3=\vtop{{\obeylines\hsize=\xs#4\par}}
     \setbox\TheAdd=\vtop{\hbox to\hsize{\hfill\copy0\hfill\copy1\hfill}
                \vskip20pt\hbox to\hsize{\hfill\copy2\hfill\copy3\hfill}}
\fi\fi\fi\catcode'015=5}}\gdef\Address{\obeylines\GetAddress}

\begin{document}
%%%%% ------------- fill in your data below this line  -------------------
%%%%%    The following lines \Title ... \EndAddress must ALL be present
%%%%%    and in the given order.
\Title

  On the leading terms of  
  Zeta isomorphisms  and   $p$-adic $L$-functions
 in non-commutative Iwasawa theory

\ShortTitle

non-commutative $p$-adic $L$-functions

\SubTitle    
Dedicated to John Coates  
%%%%%    A possible subtitle goes here.

\Author David Burns and Otmar Venjakob

\ShortAuthor
%%%%%%   Running title for even numbered pages, ONE line, please.
%%%%%%   If none is given, \Author will be used instead.
\EndTitle

\Abstract We discuss the formalism of Iwasawa theory descent in the setting of the
localized $K_1$-groups of Fukaya and Kato. We then prove interpolation formulas
for the `leading terms' of the global Zeta isomorphisms that are associated to
certain Tate motives and of the $p$-adic $L$-functions that are associated to
certain critical motives. 
 \EndAbstract

\MSC Primary 11G40; %(L-function of varieties);
Secondary
11R65 %Picard and Class of Orders
19A31 %K0 of orders
19B28 %K1 of orders.
\EndMSC
\KEY
%%%%%    Keywords and Phrases:
\EndKEY
%%%%%    All 4 \Address lines below must be present. To center the last
%%%%%    entry, no empty lines must be between the following \Address
%%%%%    and \EndAddress lines.

\Address David Burns
King's College London,
Dept. of Mathematics,

London WC2R 2LS,
United Kingdom.

\Address Otmar Venjakob
Universit\"{a}t Bonn
Mathematisches Institut
Beringstra{\ss}e 1
D-53115 Bonn
Germany.
\Address
\Address
\EndAddress

\section{Introduction}

In the last few years there have been several significant developments in
non-commutative Iwasawa theory. 

Firstly, in \cite{cfksv}, Coates, Fukaya, Kato, Sujatha and the second named author
formulated a main conjecture for elliptic curves without complex
multiplication. More precisely, if $F_\infty$ is any Galois extension of a number
field $F$ which contains the cyclotomic $\bz_p$-extension $F_{\rm cyc}$ of $F$
and is such that $\Gal(F_\infty/F)$ is a compact $p$-adic Lie group with no
non-trivial $p$-torsion, then Coates et al formulated a $\Gal(F_\infty/F)$-equivariant
main conjecture for any elliptic curve which is defined over $F$,  has good
reduction at all places above $p$ and whose Selmer group (over $F_\infty$) satisfies a
certain natural torsion condition. 

Then, in \cite{fukaya-kato}, Fukaya and Kato formulated a natural main conjecture for any
compact $p$-adic Lie extension of $F$ and any critical motive $M$ which is
defined over $F$ and has good ordinary reduction at all places above $p$.   

The key feature of \cite{cfksv} is the use of the localization sequence of algebraic
$K$-theory with respect to a canonical Ore set. However, the more general
approach of \cite{fukaya-kato}  is rather more involved and uses a notion of `localized
$K_1$-groups' together with Nekov\'a\v{r}'s theory of Selmer complexes and the
(conjectural) existence of certain canonical $p$-adic $L$-functions.  See \cite{ven-BSD} for a survey.

The $p$-adic $L$-functions of Fukaya and Kato satisfy an interpolation formula
which involves both the `non-commutative Tamagawa number conjecture' (this is a
natural refinement of the `equivariant Tamagawa number conjecture' formulated
by Flach and the first named author in \cite{bufl01}
and hence also implies the `main conjecture of non-abelian Iwasawa theory'
discussed by Huber and Kings in \cite{hu-ki}) as well as a
local analogue of the non-commutative Tamagawa number conjecture.  Indeed, by
these means, at each continuous finite dimensional $p$-adic representation
$\rho$ of $\Gal(F_\infty/F)$, the `value at $\rho$' of the $p$-adic $L$-function is
explicitly related to the value at the central critical point of the complex
$L$-function associated to the `$\rho^*$-twist' $M(\rho^*)$ of $M$. However, if the
Selmer module of $M(\rho^*)$ has strictly positive rank (and by a recent result of Mazur and Rubin  \cite{mazur-rubin}, which is itself equivalent to a special case of an earlier result of Nekov\'a\v{r} \cite[Th.\ 10.7.17]{nek}, this should often be the case),
then both sides of the Fukaya-Kato interpolation formula are equal to zero. 

The main aim of the present article is therefore to extend the formalism of
Fukaya and Kato in order to obtain an interesting interpolation formula for all
representations $\rho$ as above. To this end we shall introduce a notion of
`the leading term at $\rho$' for elements of suitable localized $K_1$-groups.
This notion is defined in terms of the Bockstein homomorphisms that have
already played significant roles (either implicitly or explicitly) in work of
Perrin-Riou \cite{perrin93,perrin2000},  of Schneider \cite{schneider85,
schneider83,schneider82,schneider81} and of Greither and the first named author \cite{burns-greither2003, fgt} and have been systematically
incorporated into Nekov\'a\v{r}'s theory of Selmer complexes \cite{nek}. We then give two explicit applications of this approach in
the setting of extensions $F_\infty/F$ with $F_{\rm cyc}\subseteq F_\infty$. We show first
that the `$p$-adic Stark conjecture at $s = 1$', as formulated by Serre \cite{ser} and interpreted by Tate in \cite{tate}, can be
reinterpreted as providing interpolation formulas for the leading terms of the
global Zeta isomorphisms associated to certain Tate motives in terms of the
leading terms at $s= 1$ (in the classical sense) of  the $p$-adic Artin
$L$-functions that are constructed by combining Brauer induction with the
fundamental results of Deligne and Ribet and of Cassou-Nogues. We then also
prove an interpolation formula for the leading terms of the Fukaya-Kato
$p$-adic $L$-functions which involves the leading term at the central critical
point of the associated complex $L$-function, the Neron-Tate pairing and
Nekov\'a\v{r}'s $p$-adic height pairing.         

In a subsequent article we shall apply the approach developed here to describe
the leading terms of the `algebraic $p$-adic $L$-functions' that are introduced
by the first named author in \cite{burns-refinedMC}, and we shall use the resulting description to prove that the main
conjecture of Coates et al for an extension $F_\infty/F$ and an elliptic curve $E$
implies the equivariant Tamagawa number conjecture for the motive $h^1(E)(1)$
at each finite degree subextension of $F_\infty/F$. We note that this result provides
a partial converse to the theorem of Fukaya and Kato which shows that, under a
natural torsion hypothesis on Selmer groups, the main conjecture of Fukaya and
Kato specialises to recover the main conjecture of Coates et al.       

The main contents of this article are as follows. In \S2 we recall some basic
facts regarding (non-commutative) determinant functors and the localized
$K_1$-groups of Fukaya and Kato. In \S3 we discuss the formalism of Iwasawa
theory descent in the setting of localized $K_1$-groups and we introduce a
notion of the leading terms at $p$-adic representations for the elements of
such groups. We explain how this formalism applies in the setting of the
canonical Ore sets introduced by Coates et al, we show that it can be
interpreted as taking values after `partial derivation in the cyclotomic
direction', and we use it to extend several well known results concerning
Generalized Euler Poincare characteristics. In \S4 we recall the `global Zeta
isomorphisms' that are conjectured to exist by Fukaya and Kato, and in \S5 
we prove an interpolation formula for the leading terms of the global Zeta
isomorphisms that are associated to certain Tate motives. Finally, in \S6, we
prove an interpolation formula for the leading terms of the $p$-adic
$L$-functions that are associated to certain critical motives.  

We shall use the same notation as in \cite{ven-BSD}.\\

It is clear that the recent developments in non-commutative Iwasawa theory are
due in large part to the energy, encouragement and inspiration of John Coates.
It is therefore a particular pleasure for us to dedicate this paper to him on
the occasion of his sixtieth birthday. 

This collaboration was initiated during the conference held in Boston in June
2005 in recognition of the sixtieth birthday of Ralph Greenberg. The authors
are very grateful to the organizers of this conference for the opportunity to
attend such a stimulating meeting.

\section{Preliminaries} \label{prelim}

\subsection{Determinant functors}
For any ring $R$ we write ${  B}(R)$ for the category of bounded
complexes of (left) $R$-modules, ${  C}(R)$ for the  category of
bounded complexes of finitely generated (left) $R$-modules, $P(R)$
for the category of finitely generated projective (left)
$R$-modules, $ {C}^p(R)$ for the category of bounded
(cohomological) complexes of finitely generated projective (left)
$R$-modules. By $D^p(R)$ we denote the category of perfect
complexes as full triangulated subcategory of the derived category
$D^b(R)$ of the homotopy category of $B(R).$ We write
$(C^p(R),{\rm quasi})$ and $(D^p(R),{\rm is})$ for the subcategory
of quasi-isomorphisms of $C^p(R)$ and isomorphisms of $D^p(R),$
respectively.

For each complex $C = (C^\cdot,d_C^\cdot)$ and each integer $r$ we
define the $r$-fold shift $C[r]$ of $C$ by setting $C[r]^i=
C^{i+r}$ and $d^i_{C[r]}=(-1)^rd^{i+r}_C$ for each integer $i$.

We first recall that for any (associative unital) ring $R$ there
exists a Picard category $\C_R$ and a determinant functor
\[\d_R:( {C}^p(R),{\rm quasi})\to \C_R\]
with the following properties (for objects $C,C'$ and $C''$ of
$\mathrm{C}^p(R)$)

\begin{itemize}
\item[d)]\footnote{The listing starts with d) to be compatible with the notation of  \cite{ven-BSD} where a)-c) describe properties of the category $\C_R.$ } If $0\to C'\to C\to C''\to 0$ is a short exact sequence
of complexes, then there is a canonical isomorphism
\[\d_R(C)\cong\d_R(C')\d_R(C'')\]
which we take as an identification.
\item[e)] If $C$ is acyclic, then the quasi-isomorphism $0\to C$
induces a canonical isomorphism
\[\u_R\to\d_R(C).\]
\item[f)] For any integer $r$ one has
$\d_R(C[r])=\d_R(C)^{(-1)^r}$.
\item[g)] the functor $\d_R$ factorizes over the image of
$\mathrm{C}^p(R)$ in the category of perfect complexes
$\mathrm{D}^p(R),$ and extends (uniquely up to unique
isomorphisms)  to $(\mathrm{D}^p(R),{\rm is}).$
\item[h)] For each $C\in \mathrm{D}(R)$ we write $H(C)$ for the
 complex which has $H(C)^i = H^i(C)$ in each degree $i$ and in which all differentials are $0$. If
 $H(C)$ belongs to $\mathrm{D}^p(R)$ (in which case we shall say that $\cC$ is {\em cohomologically perfect}), then there are canonical
 isomorphisms
\[\d_R(C) \cong \d_R(H(C)) \cong \prod_{i\in \bz} \d_R(H^i(C))^{(-1)^i}.\]

\item[i)] If $R'$ is any further ring and $Y$ an $(R',R)$-bimodule
which is both finitely generated and projective as an $R'$-module,
 then the functor $Y\otimes_R-:\mathrm{P}(R)\to P(R')$ extends to a
commutative diagram
\[\begin{CD}
(\mathrm{D}^p(R), {\rm is}) @> \d_R >> \C_R\\
@V {Y\otimes_R^\mathbb{L}-}VV @VV Y\otimes_R- V\\
(\mathrm{D}^p(R'), {\rm is}) @> \d_{R'} >> \C_{R'}.
\end{CD}\]
In particular, if $R\to R'$ is a ring homomorphism and $C\in
\mathrm{D}^p(R),$ we often write $\d_R(C)_{R'}$ in place of
$R'\otimes_R\d_R(C).$
\end{itemize}

\begin{remark}\label{signs} {\em Unless $R$ is a
regular ring, property d) does not in general extend to arbitrary
distinguished triangles. The second displayed isomorphism in h) is
induced by properties d) and f). However, whilst a precise
description of the first isomorphism in h) is important for the
purposes of explicit computations, it is actually rather difficult
to find in the literature. Here we use the description given by
Knudsen in \cite[\S3]{knudsen}.} \end{remark}

%\begin{remark}\label{eltdivisor}{\em
% Let $\O=\O_L$ be the ring of integers of a finite extension $L$ of $\qp$ and let $A$ be a finite $\O$-module. For any trivialization $a:\u_\O\to \d_\O(A),$  we obtain a canonical element $c=c(a)\in L^\times=\mathrm{Aut}_{\C_L}(\u_L)$
%\[\xymatrix{
%  {\d_L} \ar[r]^{a_L} & L\otimes_\O\d_\O(A) \ar[r]^{ } & {\d_L(L\otimes_\O A)} \ar[r]^{acyc} & {\d_L}   }\] where the map 'acyc' is induced by property e). As an immediate consequence of the elementary divisor theorem one checks that $\mathrm{ord}_L(c)=\mathrm{length}_\O(A).$ This holds in particular  for a trivialization - which we denote by 'def' since it is unique up to unique isomorphism - that is induced by some exact sequence of the form $\xymatrix@C=0.5cm{
%  0 \ar[r] & {\O^n} \ar[rr]^{ } && {\O^n} \ar[rr]^{ } && A \ar[r] & 0 }$, which always exist.}
%\end{remark}

\begin{remark}{\em \label{inverse}
We have to distinguish  between at least two inverses of a map
$\phi:\d_R(C)\to \d_R(D)$ with  $C,D\in \mathrm{C}^p(R).$ The
inverse with respect to composition will be denoted by
$\overline{\phi}:\d_R(D)\to\d_R(C)$  while $\phi^{-1}:=\overline{\id_{\d_R(D)^{-1}}\cdot\phi\cdot
\id_{\d_R(C)^{-1}}}:\d_R(C)^{-1}\to \d_R(D)^{-1}$ is   the unique isomorphism such that $\phi\cdot\phi^{-1}=\id_{\u_R}$ under the identification $\d_R(X)\cdot\d_R(X)^{-1}=\u_R$ for $X=C,D.$ If  $D=C,$ then $\phi:\d_R(C)\to \d_R(C)$ corresponds uniquely to  an element of $K_1(R)\cong\mathrm{Aut}_{\C_R}(\u_R)$ by the rule $\phi\cdot \id_{\d_R(C)^{-1}}:\u_R\to \u_R.$ Under this identification  $\overline{\phi}$ and $\phi^{-1}$ agree in $K_1(R)$ and   are inverse to $\phi.$ Furthermore, the following relation between
$\circ$ and $\cdot$ is easily verified: if $\phi: A \to B$ and $\psi : B \to C$ are morphisms in $\C_R,$ then one has $\psi \circ \phi = \psi\cdot\phi\cdot {\rm id}_{B^{-1}}.$  }
\end{remark}

We shall use the following\\  

{\bf Convention:} If $\phi:\u\to A$ is a morphism and $B$ an object in $\C_R,$ then we write $\xymatrix{
  B \ar[r]^{\cdot\;\phi} & B\cdot A   }$ for the morphism $\id_B\cdot\phi.$ In particular, any morphism $\xymatrix{
  B \ar[r]^{\phi} & A   }$ can be written as $\xymatrix{
  B \ar[rr]^{\cdot\;(\id_{B^{-1}}\cdot\;\phi)} &&   A   }.$

\begin{remark}{\em \label{trivializedident}  The determinant of the complex $C=[P_0\stackrel{\phi }{\to}
P_1]$ (in degree $0$ and $1$) with $P_0=P_1=P$ is by definition
$\xymatrix@C=0.5cm{
  { \d_R(C)}\ar@{=}[r]^<(0.3){def} & {\u_R}   }$ and is defined even if $\phi$ is not an
isomorphism (in contrast to $\d_R(\phi)$). But if $\phi$ is an isomorphism, i.e.\ if $C$ is acyclic, then by e) there is
also a canonical map $\xymatrix@C=0.5cm{
  {  \u_R}\ar[r]^<(0.3){acyc} & {\d_R(C)}   },$ which is in fact  nothing else then
\[\xymatrix@C=0.5cm{ {\u_R}\ar@{=}[r] & {\d_R(P_1)\d_R(P_1)^{-1}
}\ar[rrr]^{\d(\phi)^{-1}\cdot \id_{\d(P_1)^{-1}}} &&&
{\d_R(P_0)\d_R(P_1)^{-1}}   \ar@{=}[r] & {\d_R(C)} }\] (and  which
depends in contrast to the first identification on $\phi$). Hence, by Remark \ref{inverse}
the composite $\xymatrix@C=0.5cm{
  {  \u_R}\ar[r]^<(0.3){acyc} & {\d_R(C)} \ar@{=}[r]^<(0.4){def} & {\u_R}
  }$ corresponds to $\d_R(\phi)^{-1}\in K_1(R).$  
 In order to distinguish between the above identifications between $\u_R$ and $\d_R(C)$     we also say that $C$ is {\em trivialized by
the identity } when we refer to $\xymatrix@C=0.5cm{
  { \d_R(C)}\ar@{=}[r]^<(0.4){def} & {\u_R}   }$ (or its inverse with
  respect to composition).  Obviously,  if $\phi=\id_P,$ then the above  identifications
    coincide.}
 \end{remark}
 
\begin{remark}\label{eltdivisor} {\em 
 Let $\O=\O_L$ be the ring of integers of a finite extension $L$ of $\qp$ and let $A$ be a finite $\O$-module. For any morphism of the form  $a:\u_\O\to \d_\O(A),$ in particular  for that induced by some exact sequence of the form $\xymatrix@C=0.5cm{
  0 \ar[r] & {\O^n} \ar[r]^{ } & {\O^n} \ar[r]^{ } & A \ar[r] & 0 }$, we obtain a canonical element $c=c(a)\in L^\times=\mathrm{Aut}_{\C_L}(\u_L)$ by means of the
composite \[\xymatrix{
  {\u_L} \ar[r]^<(0.2){a_L} & L\otimes_\O\d_\O(A) \ar@^{=}[r]^{ } & {\d_L(L\otimes_\O A)} \ar[r]^<(0.3){acyc} & {\u_L}   }\] where the map 'acyc' is induced by property e). As an immediate consequence of the elementary divisor theorem one checks that $\mathrm{ord}_L(c)=\mathrm{length}_\O(A).$}
\end{remark}

\subsection{The localized $K_1$-group}
In \cite{fukaya-kato} a {\em localized} $K_1$-group was defined
for any full subcategory $\Sigma$ of $C^p(R)$ which satisfies the
following four conditions:

\begin{itemize}

\item[(i)] $0\in \Sigma,$

\item[(ii)] if $C,C'$ are in $C^p(R)$ and $C$ is quasi-isomorphic
to
$C',$ then $C\in \Sigma$ $\Leftrightarrow$ $C' \in \Sigma,$

\item[(iii)] if $C\in \Sigma,$ then also $C[n]\in\Sigma$ for all
$n\in \bz,$

\item[(iv)] if $0 \to C' \to C \to C'' \to 0$ is an exact sequence
in $C^p(R)$ with $C',C''\in\Sigma$ then also
  $C\in \Sigma.$

\end{itemize}

Since we want to apply the same construction to a subcategory
which is not necessarily closed under extensions, we weaken the
last condition to

\begin{itemize}
\item[(iv$'$)]  if  $C'$ and $C''$ belong to  $\Sigma$, then $
C'\oplus C''$ belongs to $\Sigma$.
\end{itemize}

\begin{defn}{\em (Fukaya-Kato) Assume that $\Sigma$ satisfies (i), (ii), (iii) and (iv$'$). The {\em localized $K_1$-group} $K_1(R,\Sigma)$ is defined to be the (multiplicatively written)
abelian group which has as generators symbols of the form $[C,a]$
for each $C\in \Sigma$ and morphism $a:\u_R\to \d_R(C)$ in $\C_R$
and relations
\begin{itemize}
\item[(0)] $[0,\id_{\u_R}]=1,$
\item[(1)] $[C',\d_R(f)\circ a]=[C,a]$ if $f:C\to C'$ is an
quasi-isomorphism with $C$ (and thus $C'$) in $\Sigma,$
\item[(2)] if $0 \to C' \to C \to C'' \to 0$ is an exact sequence
in $\Sigma,$ then
\[ [C,a]=[C',a']\cdot [C'',a'']\]
where $a$ is the composite of $a'\cdot a''$ with the isomorphism
induced by property d),
\item[(3)] $[C[1],a^{-1}]=[C,a]^{-1}.$
\end{itemize}}
\end{defn}

\begin{remark} {\em Relation (3) is a simple consequence
of the relations (0)-(2). Note also that this definition of
$K_1(R,\Sigma)$ makes no use of the conditions (iii) and (iv$'$)
that the category $\Sigma$ is assumed to satisfy. In particular,
if $\Sigma$ satisfies (iv) (rather than only (iv$'$)), then the
above definition coincides with that given in
\cite[\S1]{fukaya-kato}. We shall often refer to a morphism in
$\C_R$ of the form $a:\u_R\to \d_R(C)$ or $a: \d_R(C)\to\u_R$ as a {\em trivialization}
(of $C$).}\end{remark}

We now assume given a left denominator set $S$ of $R$ and we let
$R_S:= S^{-1}R$ denote the corresponding localization and
$\Sigma_S$ the full subcategory of $C^p(R)$ consisting of all
complexes $C$ such that $R_S\otimes_R C$ is acyclic. For any $C\in
\Sigma_S$ and any morphism $a: \u_{R} \to \d_{R}(C)$ in $\C_R$ we
write $\theta_{C,a}$ for the element of $K_1(R_S)$ which
corresponds under the canonical isomorphism $K_1(R_S) \cong {\rm
Aut}_{\C_{R_S}}(\u_{R_S})$ to the composite
\[ \u_{R_S} \xrightarrow{} \d_{R_S}(R_S\otimes_R C)
\rightarrow \u_{R_S}\]
where the first arrow is induced by $a$ and the second by the
 fact that $R_S\otimes_RC$ is acyclic. Then it can be shown that the assignment $[C,a] \mapsto
\theta_{C,a}$ induces an isomorphism of groups
\[ {\rm ch}_{R,\Sigma_S}: K_1(R,\Sigma_S) \cong K_1(R_S)\]
(cf. \cite[prop.\ 1.3.7]{fukaya-kato}). Hence, if $\Sigma$ is any
subcategory of $\Sigma_S$ we also obtain a composite homomorphism
\[ {\rm ch}_{R,\Sigma}: K_1(R,\Sigma) \to K_1(R,\Sigma_S) \cong K_1(R_S).\]

In particular, we shall often use this construction in the following
case: $\cA \in \Sigma_S$ and $\Sigma$ is equal to smallest full subcategory
$\Sigma_{\cA}$ of  $C^p(R)$ that contains $\cA$ and also satisfies the conditions
$(i),$ $ (ii),$ $ (iii)$ and $(iv)$ that are described above. (With this definition, it
is easily seen that $\Sigma_{\cA} \subset \Sigma_{\rm ss}$).
 
%In the case that $S$ is equal to the full set of non-zero divisors
%in $R$ we abbreviate ${\rm ch}_{R,\Sigma_S}$ to ${\rm ch}_R$.

\section{Leading terms}\label{bock}

In this section we define a notion of the leading term at a continuous finite
dimensional p-adic representation of elements of suitable localized
$K_1$-groups. To do this we introduce an appropriate `semisimplicity'
hypothesis and use a natural construction of Bockstein homomorphisms. We also
discuss several alternative characterizations of this notion. We explain how
this formalism applies in the context of the canonical localizations introduced
in \cite{cfksv} and we use it to extend several well known  results
concerning Generalized Euler Poincare characteristics.

\subsection{Bockstein homomorphisms} Let $G$ be a compact $p$-adic Lie group which contains a closed normal subgroup
 $H$ such that the quotient group $\Gamma:=G/H$ is topologically isomorphic to
 $\zp.$ Then we fix a topological generator $\gamma$ of $\Gamma$ and denote by
\[ \theta\in \H^1(G,\zp)=\Hom_{\rm cont}(G,\zp)\]
the unique homomorphism $G\twoheadrightarrow \Gamma \to \zp$ which
sends $\gamma$ to $1.$ We write $\La(G)$ for the Iwasawa algebra
of $G.$ Since $\H^1(G,\zp)\cong \mathrm{Ext}^1_{\La(G)}(\zp,\zp)$
by \cite[prop.\ 5.2.14]{nsw}, the element $\theta$ corresponds to
a canonical extension of $\La(G)$-modules
\be \label{theta} 0 \to \zp \to E_\theta \to \zp \to 0. \ee
Indeed, one has $E_\theta= \bz_p^2$ upon which $G$ acts by
$\begin{pmatrix}
  1  & \theta  \\
  0 & 1
\end{pmatrix}.$

For any $\cA$ in $B(\La(G))$ we endow the complex
 $\cA\otimes_{\zp} E_\theta$ with the natural diagonal $G$-action. Then (\ref{theta})
 induces an exact sequence in $B(\La(G))$ of the form
\[ 0 \to \cA \to \cA\otimes_{\bz_p} E_\theta \to \cA \to 0. \]
In $D^b(\La(G))$ we thus obtain a morphism, the cup product by
$\theta$ (depending on the choice of $\gamma$),
\[ \cA \xrightarrow{\theta\cup-} \cA [1] \]
which we also denote just by $\theta.$

Now let $\rho: G\to GL_n(\O)$ be a (continuous) representation of
$G$ on $T_\rho=\O^n,$ where $\O=\O_L$ denotes the ring of integers
of a finite extension $L$ of $\qp.$ We will be mainly interested
in the morphism induced by $\theta$ after forming the derived
tensor product with $\O^n$ considered as right $\La(G)$-module via
the transpose $\rho^t$ of $\rho$
\[ \O^n \otimes_{\La(G)}^\mathbb{L}\cA \xrightarrow{\theta_*} \O^n \otimes_{\La(G)}^\mathbb{L}\cA
[1].\]
In particular, we have homomorphisms
\[ \mathbb{T}\mathrm{or}_i^{\La(G)}(T_\rho,\cA) \xrightarrow{H^{-i}(\theta_*)}
\mathbb{T}\mathrm{or}_{i-1}^{\La(G)}(T_\rho,\cA)\]
of the hypertor groups
$\mathbb{T}\mathrm{or}_i^{\La(G)}(T_\rho,\cA):=H^{-i}(\O^n
\otimes_{\La(G)}^\mathbb{L}\cA). $ We refer to the map
\[ \B_i:=H^{-i}(\theta_*)\]
as the {\em Bockstein morphism (in degree $i$)}. %and we note that
%this morphism is independent of the choice of $\gamma.$

\subsection{The case $G=\Gamma$}\label{g=gamma}

In this section we consider the case $G=\Gamma$ and take the
trivial $\Gamma$-module $\zp$ for $\rho.$ We set $T:=\gamma-1.$

For any complex $\cA\in B(\La(\Gamma))$ it is clear that the
canonical short exact sequence
\[ 0
\to \La(\Gamma) \xrightarrow{\times T}\La(\Gamma) \to \bz_p\to 0\]
induces an exact triangle in $D^p(\La(\Gamma))$ of the form
\be \label{T} \cA \xrightarrow{\times T} {\cA} \to
 {\zp}\otimes^\mathbb{L}_{\La(\Gamma)}\cA \to \cA [1]. \ee
However, in order to be as concrete as possible, we choose to
describe this result on the level of complexes. To this end we fix
the following definition of the mapping cone of a map $f:\cA\to
\cB$ of complexes:
\[ \cone(f):=\cB \oplus \cA[1],\]
with differential in degree $i$ equal to
\[  d^i_{\cone(f)}:=\begin{pmatrix}
   d^i_{\cB} & f^i \\
  0 & -d^{i+1}_{\cA}
\end{pmatrix}: B^i\oplus A^{i+1}\to B^{i+1}\oplus A^{i+2}.\]

If $\cA$ belongs to $C^p(\La(\Gamma))$, then we set
\[ \cone(\cA):=\cone({\cA} \xrightarrow{T} {\cA})\]
and
\[ {A_0^\cdot}:=\zp\otimes_{\La(G)}^\mathbb{L}\cA.\]
Then we have the following morphism of complexes $\pi : \cone(\cA)
\to A_0^\cdot$
\[\begin{CD} @> >>  A^{i-1}\oplus A^i @> d_\cone^{i-1}>> A^{i}\oplus A^{i+1} @> d_\cone^{i}>>
A^{i+1}\oplus A^{i+2} @> d_\cone^{i+1}>> \\
@. @V \pi^{i-1} VV @V \pi^{i}VV @V \pi^{i+1}VV \\
 @> >> {A_0^{i-1}} @> d^{i-1}_{A_0^\cdot}>> A_0^{i} @>
d^{i}_{A_0^\cdot}>> {A_0}^{i+1} @> d^{i+1}_{A_0^\cdot}
>>
\end{CD}\]
where, in each degree $i$, $\pi^i$ sends $(a,b) \in A^i \oplus
A^{i+1}$ to the image of $a$ in $A^i/TA^i \cong A_0^i$. It is easy
 to check that $\pi$ is a quasi-isomorphism.

Now from \eqref{T} we obtain short exact sequences
\be \label{bock-decomp} 0 \to \H^i(\cA)_\Gamma \to \mathbb{H}_{-i}(\Gamma,\cA) \to \H^{i+1}(\cA)^\Gamma \to 0 \ee
where ${\mathbb{H}_{i}(\Gamma,\cA)}:=
\mathbb{T}\mathrm{or}_i^{\La(\Gamma)}(\zp,\cA)$ denotes hyper
group homology of $\cA$ with respect to $\Gamma.$ Furthermore, for
any $\La(\Gamma)$-module $M$ we write $M_\Gamma=M/TM$ and
$M^\Gamma={_T}M$ (kernel of multiplication by $T$) for the maximal
quotient module, resp. submodule, of $M$ upon which $\Gamma$ acts
trivially.

\begin{lem}\label{bockstein}
In each degree $i$ the homomorphism
$\B_i:{\mathbb{H}_{i}(\Gamma,\cA)}\to
{\mathbb{H}_{i-1}(\Gamma,\cA)}$ can be computed as the composite
\[ {\mathbb{H}_{i}(\Gamma,\cA)} \to {\H^{-i+1}(\cA)^\Gamma} \xrightarrow{\kappa^{-i+1}(\cA)} {\H^{-i+1}(\cA)_\Gamma}
 \to {\mathbb{H}_{i-1}(\Gamma,\cA)}\]
where the map
\[ \kappa^i :{\H^{i}(\cA)^\Gamma} \hookrightarrow
 \H^{i}(\cA)\twoheadrightarrow{\H^{i}(\cA)_\Gamma}\]
is induced by the identity.
\end{lem}

\begin{proof} As shown in \cite[lem.\ 1.2]{Ra-Zi}, on the level of complexes
$\theta$ is given by the map
\[ \theta:\cone(\cA)\to \cone(\cA)[1]\]
which sends $(a,b)\in A^i\oplus A^{i+1}$ to $(b,0)\in
A^{i+1}\oplus A^{i+2}.$ Now let $\bar{a}$ be in $\ker
d_{A_0^\cdot}^{-i}$ representing a class in
${\mathbb{H}_{i}(\Gamma,\cA)}.$ Then there exists $(a,b)\in \ker
(d^{-i}_\cone)$ with $\pi^{-i}((a,b))=\bar{a}.$ Since $(a,b)\in
\ker (d^{-i}_\cone)$ one has $b\in \ker (d_\cA^{i+1})$ and
$Tb=-d^i_\cA(a)$. This implies that $d^i_\cA(a)$ is divisible by
$T$ (in $A^{i+1}$) and that $b=-T^{-1}d^i_\cA(a)\in A^{i+1}$. Thus
$\theta$ maps $(a,b)$ to $(-T^{-1}d^i_\cA(a),0)$ and the class in
${\mathbb{H}_{i-1}(\Gamma,\cA)}$ is represented by
$-\overline{T^{-1}d^i_\cA(a)}\in\ker (d_{A_0^\cdot}^{-i+1}).$ By
using the canonical exact sequence
\[ 0 \to {\cA}\to \cone(\cA) \to \cA [1] \to 0 \]
one immediately verifies that $\B_i$ coincides with the map
described in the lemma.
\end{proof}

From this description it follows immediately that the pair
\be \label{bockcomplex}( {\mathbb{H}_{i}(\Gamma,\cA)},\B_i)\ee
forms a (homological) complex, which by reindexing we consider as
cohomological complex if necessary. We refer to the morphism $\B_i$ defined here as the Bockstein morphism in degree $i$ of the pair $(\cA, T)$.

\begin{defn} (Semisimplicity)\label{defsemi} {\em We shall say that a complex $\cA\in C^p(\La(\Gamma))$ is {\em semisimple} if
 the cohomology of the associated complex \eqref{bockcomplex} is $\bz_p$-torsion (and hence finite) in all
 degrees. We let $\Sigma_{\rm ss}$ denote the full subcategory of
$C^p(\La(\Gamma))$ consisting of those complexes that are
semisimple. For any $\cA\in \Sigma_{\rm ss}$ we set}
\[ r_\Gamma(\cA) := \sum_{i\in \bz} (-1)^{i+1} \dim_{\qp}
(\H^i(\cA)^\Gamma \otimes_{\bz_p}{\qp}) \in \bz.\]
\end{defn}

\begin{remark}\label{ss-remark}{\em (i) If $\cA\in \Sigma_{\rm ss}$, then
  the cohomology of $\cA$ is a torsion $\Lambda(\Gamma)$-module
 in all degrees.

(ii) In each degree $i$ Lemma \ref{bockstein} gives a canonical
exact sequence
\[ 0 \to \cok(\kappa^{-i}_{\phantom{i}}) \to \ker(\B_{i})/\mathrm{im}(\B_{i+1}) \to
\ker(\kappa^{-i+1}_{\phantom{i}})\to 0.\]
This implies that a complex $\cA\in C^p(\La(\Gamma))$ belongs to
$\Sigma_{\rm ss}$ if and only if the map
$\kappa^i(\cA)\otimes_{\zp}\qp$ is bijective in each degree $i$,
and hence also that in any such case one has}
\[
 r_\Gamma(\cA) = \sum_{i\in \bz} (-1)^{i+1} \dim_{\qp} (\H^i(\cA)_\Gamma
\otimes_{\zp}\qp).
\]
\end{remark}

\begin{defn} (The canonical trivialization) {\em We write $(\mathbb{H}_\cdot(\Gamma,\cA),0)$ for the
complex which has $(\mathbb{H}_\cdot(\Gamma,\cA),0)^i =
\mathbb{H}_i(\Gamma,\cA)$ in each degree $i$ and in which all
differentials are the zero map. If $\cA \in \Sigma_{\rm ss}$, then
we obtain a canonical trivialization of
$\zp\otimes_{\La(\Gamma)}^\mathbb{L}\cA$ by taking the composite
\begin{multline}\label{triv}
t(\cA):\d_{\zp}(\zp\otimes_{\La(\Gamma)}^\mathbb{L}\cA)_{\qp}\cong
\d_{\zp}((\mathbb{H}_\cdot(\Gamma,\cA),0))_{\qp}\\=\d_{\zp}((\mathbb{H}_\cdot(\Gamma,\cA),\B_\cdot))_{\qp}\cong
\u_{\qp}\end{multline}
where the first, resp. last, isomorphism uses property h), resp.
e), of the functor $\d_{\zp}$.}
\end{defn}

\begin{remark}\label{finitecase} {\em If $\qp\otimes_{\La(\Gamma)}^\mathbb{L}\cA$ is acyclic, then $t(\cA)$
coincides with the trivialization obtained by applying property e)
directly to $\qp\otimes_{\La(\Gamma)}^\mathbb{L}\cA.$} \end{remark}

The category $\Sigma_{\rm ss}$ satisfies the conditions (i), (ii),
(iii) and (iv$'$) that are described in \S\ref{prelim} (but does
not satisfy condition (iv)). In addition, as the following lemma
shows, the above constructions behave well on
 exact triangles of semisimple complexes.

\begin{lem}\label{triangle}
 Let $\cA,\cB$ and $\cC$ be objects of $\Sigma_{\rm ss}$ which together lie in an exact triangle in $D^p(\La(\Gamma))$ of the form
\[ \cA\to \cB\to \cC\to
 \cA[1].\]
Then one has
\[ r_\Gamma(\cB)=r_\Gamma(\cA)+r_\Gamma(\cC)\]
and, with respect to the canonical isomorphism
\[ \d_{\zp}(\zp\otimes_{\La(G)}^\mathbb{L} B^\cdot)_{\qp}=\d_{\zp}(\zp\otimes_{\La(G)}^\mathbb{L}A^\cdot)_{\qp}\d_{\zp}(\zp\otimes_{\La(G)}^\mathbb{L}C^\cdot)_{\qp}\]
that is induced by the given exact triangle, one has
\[ t(\cB)= t(\cA) t(\cC).\]
\end{lem}

\begin{proof} Easy and left to the reader. \end{proof}

We write $\rho_{\rm triv}$ for the trivial representation of
$\Gamma$.

\begin{defn}\label{leadingterm} (The leading term) {\em For each $\cA \in \Sigma_{\rm ss}$ and each morphism
 $a :\u_{\La(\Gamma)}\to \d_{\La(\Gamma)}(\cA)$ in
 $\C_{\La(\Gamma)}$ we define the {\em leading term $(\cA,a)^*(\rho_{\triv})$
of the pair $(\cA,a)$ at} $\rho_{\triv}$ to be equal to
$(-1)^{r_\Gamma(\cA)}$ times the element of $\bq_p\setminus \{
0\}$ which corresponds via the isomorphisms $\bq_p^\times \cong
K_1(\bq_p) \cong {\rm Aut}_{\C_{\bq_p}}(\u_{\qp})$ to the
 composite morphism
\[ \u_{\qp} \xrightarrow{\qp\otimes_{\La(\Gamma)}a} {\d_{\zp}(\zp\otimes_{\La(\Gamma)}\cA)_{\qp}} \xrightarrow{t(\cA)}
\u_{\qp}.\]
After taking into account Lemma \ref{triangle}, this construction
induces a well defined homomorphism of groups}
\begin{align*} (-)^*(\rho_{\triv}): K_1(\La(\Gamma),\Sigma_{\rm ss}) &\to
\bq_p^\times\\
[\cA,a] &\mapsto [\cA,a]^*(\rho_{\rm triv}) :=
(\cA,a)^*(\rho_{\rm triv}).\end{align*}
\end{defn}

The reason for the occurrence of $\rho_{\rm triv}$ in the above
definition will become clear in the next subsection. In the remainder
of the current section we justify the name `leading term' by
explaining the connection between $(\cA,a)^*(\rho_{\triv})$ and
the leading term of an appropriate characteristic power series.

To this end we note that $\Sigma_{\rm ss}$ is a subcategory of the
full subcategory of $C^p(\Lambda(\Gamma))$ consisting of those
complexes $C$ for which $Q(\Gamma)\otimes_{\La(\Gamma)}C$ is
acyclic (by Remark \ref{ss-remark}(i)), and hence that there
exists a homomorphism
\[ {\rm ch}_{ \Gamma }:={\rm ch}_{\La(\Gamma),\Sigma_{\rm ss}}: K_1(\La(\Gamma),\Sigma_{\rm ss}) \to
K_1(Q(\Gamma)) \cong Q(\Gamma)^\times.\]

If $L$ is any field extension of $\bq_p$, and $\O$ the valuation
ring in $L$, then the identification $\La_\O(\Gamma)\cong \O[[T]]$
(which, of course, depends on the choice of $T=\gamma-1$) allows
any element $F\in Q_\O(\Gamma)^\times$ to be written uniquely as
\be\label{decomp} F(T)=T^rG(T)\ee
with $r=r(F)\in \bz$ and $G(T)\in Q_\O(\Gamma)$ such that $G(0)\in
L^\times.$ The leading coefficient of $F$ with respect to its
 Taylor series expansion in the Laurent series ring
$L\{\{T\}\}$ is thus equal to $F^*(0):=G(0)$. Let $\p$ be the kernel of the augmentation map $\La(\Gamma)\to
\zp$ and denote by $R$ the localisation $\La(\Gamma)_\p.$ This is
a discrete valuation ring with uniformizer $T$ and the residue
class field $R/(T)$ is isomorphic to $\qp.$ The  notation of semi-simplicity extends naturally to $R$ via the exact functor $R\otimes_{\La(\Gamma)}-:$   indeed, if $\cA\in D^p(R),$ then the analogue of \eqref{T} induces a Bockstein map  for $\mathbb{T}\mathrm{or}_\cdot^R(R/(T),\cA)$ and we say
that $\cA$ is  {\em semisimple } if   the associated complex $(\mathbb{T}\mathrm{or}_\cdot^R(R/(T),\cA),\B_\cdot)$ is acyclic. 

\begin{prop}\label{commut1} Let $a: \u_{\Lambda(\Gamma)} \to
\d_{\Lambda(\Gamma)}(\cA)$  be a morphism in $\C_{\La(\Gamma)}$ with $\cA$ in $\Sigma_{\rm ss}.$  

\begin{itemize}
\item[(i)] {\em (Order of vanishing)} For $\L:=[\cA,a]$
 one has $r({\rm
ch}_{ \Gamma }(\L)) = r_\Gamma(\cA)$.
\item[(ii)] {\em (Leading terms)} One has a commutative diagram of
abelian groups
\[ \begin{CD} K_1(\La(\Gamma),\Sigma_{\rm ss}) @> {\rm ch}_{ \Gamma } >> K_1(Q(\Gamma))\\
@V {(-)^*(\rho_{\triv})}VV  @VV {(-)^*(0)}V \\
{\bq_p^\times} @= {\bq_p^\times}.\end{CD} \]
%
%we have $(-1)^{r_\Gamma(\cA)}\L^*(\rho_{tri})=F_\L^*(0).$
\end{itemize}
\end{prop}

\begin{proof}
 It is easy to see that
both
  maps $(-)^*(\rho_{\triv})$ and ${\rm ch}_{ \Gamma }$ factor via flat base change $R\otimes_{\La(\Gamma)}-$ through
  $K_1(R,\Xi),$ where $\Xi$ denotes the full subcategory of
  $C^p(R)$ consisting of those complexes which are semisimple. Thus it   suffices  to show the commutativity of the above diagram with  $K_1(\La(\Gamma),\Sigma_{\rm ss})$
  replaced by $K_1(R,\Xi).$ Moreover, by   Lemma \ref{generator} below this
  is reduced to the case where $\cA$ is a complex of the form $R \xrightarrow{a} R$ where $R$ occurs in degrees $-1$ and $0$
  and $a$ denotes multiplication by either $T$ or $1$. Further, since the complex $R \xrightarrow{\times 1} R$ is acyclic
  we shall therefore assume that $a$ denotes multiplication by $T$.

  Now $\mathrm{Mor}_{\C_R}(\u_R,\d(\cA))$ is a
    $K_1(R)$-torsor and so all possible trivialization
    arise in the following way where   $\epsilon\in R^\times:$   The module homomorphisms $R\to A_i$ which send $1$ to $1\in R=A_{-1}$ and to $\epsilon\in R=A_{0},$ respectively, induce maps
    $can_1 : \d_R(R)\to\d_R(A_{-1})$ and $can_\epsilon: \d_R(R)\to\d_R(A_0)$ which give rise to the
    trivialization $a_\epsilon:= (can_1)^{-1}\cdot can_\epsilon:
    \u_R\to \d_R(\cA).$ Setting $\L_\epsilon:= [\cA,a_\epsilon]$ one checks easily that
    $\mathrm{ch}_\Gamma(\L_\epsilon)=T^{-1}\epsilon$ and thus
    $\mathrm{ch}_\Gamma(\L_\epsilon)^*(0)=\epsilon(0).$ On the other hand, the
    bockstein map $\B_1$ is just $\zp \xrightarrow{-1} {\zp}$ as one checks using the
      description of $\B$ in the proof of Lemma \ref{bockstein}. Thus
    $\L^*_\epsilon(\rho_{\triv})$ is equal to $(-1)^{r_\Gamma(\cA)}$ times the determinant of
\[ {\qp} \xrightarrow{\epsilon(0)} {\qp}  \xrightarrow{(\B_1)_{\qp}^{-1}=-1} {\qp} \xrightarrow{1}
{\qp}.\]
Hence, observing  that
      $r_\Gamma(\cA)=-1=r(\mathrm{ch}_\Gamma(\L_\epsilon)),$  we have   $\L^*_\epsilon(\rho_{\triv}) =  \epsilon(0) $
  which proves       the Proposition.
\end{proof}

\begin{lem}\label{generator}
Let $R$ be a discrete valuation ring with uniformizer $T$ and
assume that $\cA\in C^p(R)$ is semisimple.  Then $\cA$ is  isomorphic in $  C^p(R)$ to the direct sum of
finitely many  complexes of the form   $R \to R$ where the differential is equal to
multiplication by either $1$ or $T.$ 
\end{lem}

\begin{proof}
Assume that $m$ is the maximal degree such that $A^m\neq 0$ and
fix an isomorphism $D:R^d\cong A^m.$ Let $(e_1,\ldots, e_d)$ be
the standard basis of $R^d.$ The semisimplicity of $\cA$ is equivalent to  the property that   $T\cdot\H^i(\cA)=0$
for all $i\in \bz$ as can be seen   easily by using an analogue  of   Lemma \ref{bockstein}; indeed, both conditions are equivalent to  the property that the analogues of the maps $\kappa^i$ in that Lemma are isomorphisms for all $i.$
Thus, for $1\leq i\leq d,$  $Te_i$ is in the image of $D^{-1}\circ d^{m-1}.$ We set
$h_i=1$ if $e_i$ is already in this image and $h_i=T,$ otherwise.
By $H$ we denote the diagonal matrix with the elements
$h_1,\ldots,h_d.$ Since the image of the map $R^d \xrightarrow{H}
 R^d$ coincides by definition with that of $D^{-1}\circ
  d^{m-1}$ we obtain a retraction $E:\, R^d \xrightarrow{} A^{m-1} $ (i.e.\ with left inverse "$H^{-1}\circ D^{-1}\circ d^{m-1}$") making the following
    diagram commutative
\[ \begin{CD} @> >> 0 @> >> R^d @> H >> R^d @> >> 0 @> >> \\
@. @V VV @V E VV @V D VV @V VV \\
@> >> A^{m-2} @> d^{m-2}>> A^{m-1} @> d^{m-1}>> A^m @> d^m >> 0 @>
>>  .\end{CD}\]
If $\cB$ denotes the upper row of this diagram and $\cC:=\cA/\cB$
the associated quotient complex (not the mapping cone!), one
checks readily that there exists a {\em split} exact sequence
 $\,0 \to { \cB} \to {\cA} \to {\cC} \to 0.$ Since $\cC$ is again semisimple and has a strictly shorter
  length than $\cA$ the proof is accomplished by induction.
\end{proof}

\begin{remark}{\em It will be clear to the reader that analogous statements hold for all results of this
subsection if we replace $\zp$ by $\O,$ $\qp$ by $L,$
$\La(\Gamma)$ by $\La_\O(\Gamma):=\O [[\Gamma]],$ and $\Q(\Gamma)$
by $Q_\O(\Gamma),$ the quotient field of
$\La_\O(\Gamma).$}\end{remark}

\subsection{The general case}\label{tgc}
For any continuous representation $\rho: G \to {\rm GL}_n(\O)$ we
regard the complex
\[ \cA(\rho^*):=\cA\otimes_{\zp} \O^n\]
as a complex of $\La(G)$-modules by means of the following
$G$-action: $g(a\otimes o):=ga\otimes\rho^*(g)o$ for $g\in G,$
$a\in A^i$ and $o\in \O^n.$

\begin{defn} (Semisimplicity at $\rho$) \label{defsemirho}{\em We shall say that a complex $\cA\in D^p(\La(G))$ is {\em semisimple at $\rho$}
if the cohomology of the associated complex
 $(\mathbb{H}_\cdot(G,\cA(\rho^*)),\B_\cdot)$ is $\bz_p$-torsion in each degree. We let $\Sigma_{{\rm ss}-\rho}$ denote the full subcategory
 of $C^p(\La (G))$ consisting of all complexes that
are semisimple at $\rho$, and we note that $\Sigma_{{\rm
ss}-\rho}$ satisfies the conditions (i), (ii), (iii) and (iv$'$)
that are described in \S\ref{prelim}. For each $\cA \in
\Sigma_{{\rm ss}-\rho}$ we set
\[ r_G(\cA)(\rho) := \sum_{i\in \bz}
(-1)^{i+1} \dim_L\big(\mathbb{H}_i(H,\cA(\rho^*))^\Gamma\otimes_\O
L\big) \in \bz,\]
where $\mathbb{H}_i(H,-):= \mathbb{T}{\rm or}_i^{\La(H)}(\O,-)$
denotes hyper $H$-group homology.}
 \end{defn}
 
\begin{defn} (Finiteness at $\rho$) {\em Similarly we shall say that a complex $\cA\in D^p(\La(G))$ is {\em finite at $\rho$}
if the cohomology groups
 $\mathbb{H}_i(G,\cA(\rho^*))$ are $\bz_p$-torsion in each degree. We let $\Sigma_{{\rm fin}-\rho}$ denote the full subcategory
 of $C^p(\La (G))$ consisting of all complexes that
are finite at $\rho$, and we note that $\Sigma_{{\rm fin}-\rho}$ satisfies the conditions (i), (ii), (iii) and (iv)
that are described in \S\ref{prelim}. In particular we have $\Sigma_{{\rm
fin}-\rho}\subseteq \Sigma_{{\rm
ss}-\rho}.$}
 \end{defn}

In the next result we consider the tensor product
$\La_\O(\Gamma)\otimes_\O\O^n$ as a natural
$(\La_\O(\Gamma),\La(G))$-bimodule where $\La_\O(\Gamma)$ acts by
 multiplication on the left and $\La(G)$ acts on the right via the rule
$(\tau\otimes o)g:= \tau \bar{g}\otimes \rho(g)^to$ for each $g\in
G$ (with image $\bar{g}$ in $\Gamma$), $o\in \O^n$ and $\tau\in
\La_\O(\Gamma).$    For a  complex $ \cA\in \Sigma_{{\rm
ss}-\rho}$, we write
 $t(\cA(\rho^*))$ for the trivialization $t(\cA_{\rho^*}):\d_\O(\O\otimes_{\Lambda_\O (\Gamma)}^{\mathbb L}{\cA_{\rho^*}})_L\to \u_L$   defined in \eqref{triv}. By the following Lemma this amounts to the same as taking the composite
\begin{multline}\label{triv
rho}
t(\cA(\rho^*)):\d_\O(\O\otimes^\mathbb{L}_{\Lambda_\O(G)}\cA(\rho^*) )_L\cong 
\d_{\O}((\mathbb{H}_\cdot(G,\cA(\rho^*)),0))_{L}\\=\d_{\O}((\mathbb{H}_\cdot(G,\cA(\rho^*)),\B_\cdot))_{L}\cong
\u_{L}\end{multline}
 where the first, resp. last, isomorphism uses property h), resp.
e), of the functor $\d_{\O}$.

\begin{lem}\label{tensor} Fix $\cA \in C^p(\La(G))$.

\begin{itemize}

\item[(i)] One has $\cA \in \Sigma_{{\rm ss}-\rho}$ if and only if
\[ {A_\rho^\cdot}:=(\La_\O(\Gamma)\otimes_\O\O^n)\otimes^\mathbb{L}_{\La(G)}\cA\cong
\La_\O(\Gamma)\otimes^\mathbb{L}_{\La_\O(G)}\cA(\rho^*) \]
belongs to $\Sigma_{\rm ss}$ (when considered as an object of
$C^p(\La_{\O}(\Gamma))$). In addition, in any such case the Bockstein map of $\mathbb{H}(G,A^\cdot(\rho^*))$ is equal to the Bockstein map of the pair $(\cA_\rho,T)$  and one has  natural
isomorphisms in $C^p(\O)$
\[ \O\otimes_{\Lambda_\O (\Gamma)}^{\mathbb L}{A_\rho^\cdot} \cong
\O^n\otimes^\mathbb{L}_{\La(G)}A^\cdot\] and in $B(\O).$
\[ \La_\O(\Gamma)\otimes^\mathbb{L}_{\La(G)}\cA\cong
\O\otimes^\mathbb{L}_{\La(H)}\cA.\]

\item[(ii)] If $\cA \in \Sigma_{{\rm ss}-\rho}$, then $\,
r_G(\cA)(\rho) = r_\Gamma({A_\rho^\cdot}).$

\item[(iii)] If $\cA, \cB$ and $\cC$ are objects of $\Sigma_{{\rm
ss}-\rho}$ which together lie in an exact triangle in
$D^p(\La(G))$ of the form
\[ A^\cdot\to B^\cdot\to C^\cdot\to
A^\cdot[1],\]
then one has
\[
r_G(\cB)(\rho)=r_G(\cA)(\rho)+r_G(\cC)(\rho)\]
and,  with respect to the canonical isomorphism
\[ \d_{\O}(\O\otimes^\mathbb{L}_{\Lambda_\O(G)}\cB(\rho^*) )_{L}=\d_{\O}(\O\otimes^\mathbb{L}_{\Lambda_\O(G)}\cA(\rho^*) )_{L}\d_{O}(\O\otimes^\mathbb{L}_{\Lambda_\O(G)}\cC(\rho^*) )_{L}\]
that is induced by the given exact triangle, one has
\[ t(B^\cdot(\rho^*))= t(A^\cdot(\rho^*))\cdot t(C^\cdot(\rho^*)).\]
\end{itemize}
\end{lem}

\begin{proof} Claim (i) follows from the fact that the specified actions induce natural isomorphisms in $D^p(\O)$ of the form %
\begin{eqnarray*}
\O^n\otimes_{\La(G)}^\mathbb{L}\cA&\cong&
\O\otimes_{\La_\O(G)}^\mathbb{L}\cA(\rho^*)\\
&\cong&\O\otimes_{\La_\O(\Gamma)}^\mathbb{L}\big( \La_\O(\Gamma)\otimes_{\La_\O(G)}^\mathbb{L}\cA(\rho^*)\big)\\
&\cong&\O\otimes_{\La_\O(\Gamma)}^\mathbb{L}\big((\La_\O(\Gamma)\otimes_\O\O^n)\otimes_{\La(G)}^\mathbb{L}\cA\big)
\end{eqnarray*}
and from the functorial construction of the Bockstein  map. 
Claim (ii) follows directly by using claim (i) to compare the
definitions of the terms $r_G(\cA)(\rho)$ and
$r_\Gamma({A_\rho^\cdot})$. To prove claim (iii)
 we observe that, by claim (i), the given triangle induces an exact triangle of
semisimple complexes in $D^p(\La(\Gamma))$ of the form
\[ A^\cdot_\rho\to B^\cdot_\rho\to
C^\cdot_\rho\to A^\cdot_\rho[1].\]
The equalities of claim (iii) thus follow from claim (ii) and the
results of Lemma \ref{triangle} as applied to the above triangle.
\end{proof}

\begin{defn}\label{leadingterm-rho}(The leading term at $\rho$) {\em For each complex $\cA \in \Sigma_{{\rm ss}-\rho}$
 and each morphism $a:\u_{\La(G)}\to
\d_{\La(G)}(\cA)$ we define
 the {\em leading term $(\cA,a)^*(\rho)$ of the pair $(\cA,a)$ at
$\rho$} to be equal to $(-1)^{r_G(\cA)(\rho)}$ times the element
of $L\setminus \{0\}$ which corresponds via the isomorphisms
$L^\times \cong K_1(L) \cong {\rm Aut}_{\C_L}({\u}_L)$ to the
 composite morphism
\[ {\u_L} \xrightarrow{L^n\otimes_{\La (G)}a} \d_L(L^n\otimes_{\La (G)}^\mathbb{L}\cA)
\xrightarrow{t(\cA(\rho^*))} \u_L .\]
Since
 $\Sigma_\cA \subset \Sigma_{{\rm ss}-\rho}$   Lemma
\ref{tensor} shows that the above construction induces a
well-defined homomorphism of groups}
\begin{align*} (-)^*(\rho): K_1(\Lambda(G),\Sigma_\cA) &\to L^\times\\
 [\cA,a] &\mapsto [\cA,a]^*(\rho) := (\cA,a)^*(\rho).\end{align*}
\end{defn}
If $\cA$ is clear from context, then we often write $a^*(\rho)$ in place of $[\cA,a]^*(\rho).$
It is easily checked that (in the case $G = \Gamma$ and $\rho =
\rho_{\triv}$) these definitions are compatible with those given
in \S\ref{g=gamma}. Further, in \S\ref{partial-deriv} we shall
reinterpret the expression $[\cA,a]^*(\rho)$ defined above as the
leading term at $s =0$ of a natural $p$-adic meromorphic function.

\begin{remark}\label{value}
{\em  If $\cA\in\Sigma_{\mathrm{fin}-\rho},$ then we set $ [\cA,a](\rho):=[\cA,a]^*(\rho)$ and call this the {\em value} of $[\cA,a]$ at $\rho.$ Taking into account  Remark \ref{finitecase} it is clear that this definition coincides with that given in 
 \cite[4.1.5]{fukaya-kato}.
}
\end{remark}

\subsection{Canonical localizations} We apply the constructions of
\S\ref{tgc} in the setting of the canonical localizations of $\La
(G)$ that were introduced in \cite{cfksv}.

\subsubsection{The canonical Ore sets} We recall from \cite[\S2-\S3]{cfksv} that
there are canonical left and right denominator sets $S$ and $S^*$
of
 $\La (G)$ where
\[S := \{\lambda\in
\La (G) : \La(G)/\La(G)\cdot\lambda \mbox{ is a finitely generated
} \La (H)\mbox{-module}\}\]
and
\[S^* := \bigcup_{i\geq 0} p^iS. \]
We write   $S^*$-tor for the category  of finitely
generated $\La (G)$-modules $M$ which satisfy
 $\La(G)_{S^*}\otimes_{\La
(G)}M = 0.$   We further recall from loc.\ cit.\ that a
 finitely generated $\La(G)$-module $M$ belongs to   $S^*$-tor, if and
only if  $M/M(p) $  is finitely generated when considered as a $\La
(H)$-module (by restriction) 
where $M(p)$ is the submodule of $M$ consisting of those elements
which are annihilated by some power of $p$.

%\begin{example}{\em In the setting of \S\ref{tif} it is conjectured that $SC_U(\hat{\T},\T)$ always belongs to
%$\Sigma_{S^*}$ (cf. \cite[conj.\ 5.1]{cfksv} and \cite[4.3.5 and
%prop.\ 4.3.7]{fukaya-kato}).
% Further, in the setting of \S\ref{Q(1)} the containment $A^\cdot_{K/k,T}\in \Sigma_{S}$ is in many cases equivalent to a well known conjecture of Iwasawa.
% To be more specific, we note that $A^\cdot_{K/k,T}\in \Sigma_S$ if
% and only if $\mathrm{G}(M_T(K)/K) \in S\mbox{-tor}$. In addition, if $H$ contains an open pro-$p$ subgroup $P$ which is
% both normal in $G$ and such that $G/P$ decomposes as a direct product $H/P\times \Gamma$, then Nakayama's lemma allows
% one to show that $\mathrm{G}(M_T(K)/K) \in S\mbox{-tor}$ if and only if
% the $p$-adic $\mu$-invariant of the field $(K^P)^\Gamma$ vanishes (as has been conjectured to be true by Iwasawa).}\end{example}

\subsubsection{Leading terms} We use the notation of Definition \ref{leadingterm-rho}. If $\rho: G \to {\rm
GL}_n(\O)$ is any continuous representation and $\cA$ any object
of $\Sigma_{S^*}$, then $\Sigma_\cA \subset \Sigma_{S^*}$ and so
there exists a canonical homomorphism
\[ {\rm ch}_{G,\cA}:={\rm ch}_{\La (G),\Sigma_\cA} : K_1(\La (G),\Sigma_\cA)\to
 K_1(\La(G)_{S^*}).\]
In addition, the ring homomorphism $\La(G)_{S^*}\to
M_n(Q(\Gamma))$ which sends each element $g \in G$ to
 $\rho(g)\bar{g}$ where $\bar{g}$ is the image in $\Gamma,$ induces
a homomorphism
\[ \rho_*: K_1(\La (G)_{S^*})\to K_1(M_n(Q_\O(\Gamma)))\cong
K_1(Q_\O(\Gamma)) \cong Q_\O(\Gamma)^\times.\]

\begin{prop}\label{commut2diagram} Let $\cA$ be a complex which belongs to $\Sigma_{S^*}\cap \Sigma_{{\rm ss}-\rho}$.

\begin{itemize}
\item[(i)] {\em (Order of vanishing)} One has $r_G(\cA)(\rho) =
r_\Gamma(A_\rho^\cdot)=r(\rho_*\circ {\rm ch}_{ G, \cA}(\cA))$.
\item[(ii)] {\em (Leading terms)} The following diagram of abelian
groups commutes
\[ \begin{CD} K_1(\La(G),\Sigma_\cA) @> {{\rm ch}_{ G, \cA}}>> K_1(\La(G)_{S^*})\\
@V {(-)^*(\rho)}VV @VV {(\rho_*(-))^*(0) }V \\
L^\times @= {L^\times},\end{CD} \]
where $(-)^*(0)$ denotes the `leading term' homomorphism
$K_1(Q_\O(\Gamma)) \to L^\times$ according to Proposition \ref{commut1}.
\end{itemize}
\end{prop}

\begin{proof} By Lemma \ref{tensor}(i)   we have $  \mathbb{H}_i(H,  \cA(\rho^*))=\H^{-i} ( \O\otimes^\mathbb{L}_{\Lambda(H)} \cA(\rho^*))=\H^{-i}(A_\rho^\cdot).$ Taking also into account Proposition \ref{commut1}, (i) follows from   Definition \ref{defsemi} and \ref{defsemirho}.
Claim (ii) is proved by the same arguments as used in \cite[lem.\
4.3.10]{fukaya-kato}. Indeed, one need only observe that the above
diagram arises as the following composite commutative diagram
\[ \begin{CD}
  K_1(\La(G),\Sigma_\cA) @> {{\rm ch}_{\La (G),\Sigma_\cA}}>> K_1(\La(G)_{S^*})\\
  @V {(\La_\O(\Gamma)\otimes_\O\O^n)\otimes_{\La(G)}-}VV   @VV \rho_*V \\
K_1(\La_\O(\Gamma),\Sigma_{\rm ss}) @> {\rm
ch}_{\La_\O(\Gamma),\Sigma_{\rm ss}}
>>
{K_1(Q_\O(\Gamma))}\\
@V{(-)^*(\rho_{\triv})}VV @VV{(-)^*(0)}V\\
L^\times @= L^\times \end{CD}\]
where the lower square is as in Proposition \ref{commut1}.
\end{proof}

Thus for any $F\in K_1(\La(G)_{S^*})$ we write also $F^*(\rho)$ for the {\em leading term $\rho_*(F)(0)$ of $F$ at $\rho.$}  We note that, by Proposition \ref{commut2diagram} , this notation is consistent with
that of Definition \ref{leadingterm-rho} in the case that $F$ belongs to   the image of $\mathrm{ch}_{G,\cA}.$ Similarly, we shall use the notation $F(\rho):=F^*(\rho)$ if $r(\rho_*(F))=0.$ 

\subsubsection{Partial derivatives}\label{partial-deriv} We now observe that
the constructions of the previous section allow an interpretation
of the expression $(\cA,a)^*(\rho)$ defined in \S\ref{tgc} as the
leading term at $s =0$ of a natural $p$-adic meromorphic function.

At the outset we fix a representation of $G$ of the form $\chi: G
\twoheadrightarrow \Gamma \to \bz_p^\times$ of infinite order and we set \[ c_{\chi,\gamma} :=
\log_p(\chi(\gamma)) \in \bq_p^\times.\]

We also fix $\cA \in \Sigma_{S^*}$ and a morphism $a: \u_{\La (G)}
\to \d_{\La (G)}(\cA)$, put $\L:=[\cA,a],$ and for any continuous representation
$\rho : G \to {\rm GL}_n(\O)$ we set
\[ f_{\rho}(T) := \rho_*({\rm ch}_{G,\cA}(\L)) \in K_1(Q_\O(\Gamma)) \cong
Q_\O(\Gamma)^\times.\]
Then, since the zeros and poles of elements in $Q_\O(\Gamma)$ are
discrete, the function
\[ s \mapsto f_{\L}(\rho\chi^s) := f_{\rho}(\chi(\gamma)^{s}-1)\]
is a $p$-adic meromorphic function on $\bz_p$.

\begin{lem}\label{calculus} Let $\cA$ and $a$ be as above and set $r:=r_G(\cA)(\rho).$ Then, 

\begin{itemize}
 \item[($i$)] in any sufficiently small neighbourhood of $0$
in $\bz_p$ one has
\[ \L^*(\rho\chi^s) =\L(\rho\chi^s)= f_{\L}(\rho\chi^s)\]
for all $s \not= 0,$ 
\item[($ii$)] 
    $c_{\chi,\gamma}^r \L^*(\rho)$ is the leading coefficient at $s=0$ of $f_{\L}(\rho\chi^s),$    and
\item[($iii$)]   if $r \ge 0,$
  one has
\[c_{\chi,\gamma}^r \L^*(\rho) = \frac{1}{r!}\frac{{\rm d}^r}{{\rm
d}s^r}f_{\L}(\rho\chi^s)\big|_{s=0}.\]\end{itemize}
\end{lem}

\begin{proof} In any sufficiently small neighbourhood of $0$ in $\bz_p$ one has
$f_{\rho\chi^s}(0)\in L^\times$ for all $s \not= 0$. Since
$f_{\rho\chi^s}(T) = f_{\rho}(\chi(\gamma)^{s}(T+1)-1)$ we may
therefore deduce from Proposition \ref{commut2diagram} that
$\L^*(\rho\chi^s)  = f_{\rho\chi^s}(0) =
f_{\rho}(\chi(\gamma)^{s}-1) = f_{\L}(\rho\chi^s)$ for any
such value of $s$.

In addition, if $r \ge 0$ and we factorize $f_{\rho}(T)$ as
$T^rG_{\rho}(T)$ with $G_{\rho}(T) \in Q_\O(\Gamma)$, then
$G_{\rho}(0) = f_{\rho}^*(0)$ and
\begin{eqnarray*}
\frac{1}{r!}\frac{{\rm d}^r}{{\rm
d}s^r}f_{\L}(\rho\chi^s)\big|_{s=0} &=& \lim_{0\neq s\to
0}\frac{f_{\rho}(\chi(\gamma)^{s}-1)}{s^r}\\
&=& \lim_{0\neq s\to 0}\big(\frac{(\chi(\gamma)^{s}-1)^r}{s^r}\;G_{\rho}(\chi(\gamma)^{s}-1)\big)\\
&=& \big(\lim_{0\neq s\rightarrow 0}\frac{\chi(\gamma)^{s}-1}{s}\big)^r G_{\rho}(0)\\
&=& (\log_p(\chi(\gamma)))^r f_{\rho}^*(0) \\
&=&  c_{\chi,\gamma}^r\L^*(\rho),
\end{eqnarray*}
where the last equality follows from Proposition
\ref{commut2diagram}. If $r<0$ we do not have the interpretation as partial derivative but the same arguments prove the statement concerning the leading coefficient at $s=0.$
\end{proof}

\begin{remark}
 {\em Of particular interest is the case $\chi=\chi_{\mathrm{cyc}}$ in which the above calculus can be interpreted as partial derivation in the `cyclotomic' direction (cf.\ Remark \ref{equa-lead})

.}
\end{remark}

\subsection{Generalized Euler-Poincare characteristics}
 In this section we observe that the constructions made in \S\ref{tgc} give rise to a natural extension of
 results from \cite{cfksv, fukaya-kato, ven-crelle}.

We fix a continuous representation $\rho:G \to {\rm GL}_n(\O)$
 and a complex $\cA \in \Sigma_{{\rm ss}-\rho}$ and in each degree $i$ we
set $\H^i_\B(G,-) :=
H^i\big(({\mathbb{H}_{-\cdot}(G,-)},\B_{-\cdot})\big)$. We then
 define the
 (generalized) additive, resp. multiplicative, Euler-Poincare characteristic
 of the complex $\cA(\rho^*)$ by setting
\[ \chi_{\rm add}(G,\cA(\rho^*)):=\sum_{i\in \bz} (-1)^i
\mathrm{length}_\O\big(\H^i_\B(G,\cA(\rho^*))\big),\]
resp.
\[ \chi_{\rm mult}(G,\cA(\rho^*)):=(\#\kappa_L)^{\chi_{\rm add}(G,\cA(\rho^*))}\]
where $\kappa_L$ is the residue class field of $L$. We recall that for single $\Lambda(G)$-module $M$ or rather its Pontryagin-dual $D $ similar  Euler characteristics have  also been studied by other 
authors (cf.\ \cite{cs2, zerbes, hachi-ven}). Indeed, they use the Hochschild-Serre spectral sequence to construct differentials \[d^i:\H^i(G,D)\to\H^{i}(H,D)^\Gamma\to\H^{i}(H,D)_\Gamma\to \H^{i+1}(G,D)\] where the middle map is induced by the identity; then their generalized Euler characteristic is defined similar as above using the complex $(\H^\cdot(G,D), d^\cdot)$ instead of $({\mathbb{H}_{-\cdot}(G,-)},\B_{-\cdot}).$ But by Lemma \ref{tensor} (i) one sees immediately that the Pontryagin dual of $d_i$ coincides with the Bockstein map  $\B_{i+1}:\mathbb{H}_{i+1}(G,P^\cdot)\to \mathbb{H}_i(G,P^\cdot)$ as computed in Lemma \ref{bockstein} with $\cA=\Lambda(\Gamma)\otimes^\mathbb{L}_{\Lambda(G)}P^\cdot,$ where $P^\cdot$ denotes a projective solution of $M.$ More precisely, in \cite{cs2} a truncated version of this generalized Euler characteristic is used. 
%because the conjectured vanishing of some of the higher cohomology groups in a specific arithmetic application is not known.

\begin{prop}\label{gecp} Let $\mathrm{ord}_L$ denote the valuation of $L$ which takes the value $1$ on any
 uniformizing parameter and $|-|_p$ the $p$-adic
absolute value, normalized such that $|p|_p= p^{-1}.$

If $\cA \in \Sigma_{{\rm ss}-\rho}$ and  $a:\u_{\La(G)}\to \d_{\La(G)}(\cA)$ is any morphism,
 then for $\L:=[\cA,a]$ one has
\[ \chi_{\rm add}(G,\cA(\rho^*))=\mathrm{ord}_L(\L^*(\rho))\]
and
\[ \chi_{\rm mult}(G,\cA(\rho^*))=|\L^*(\rho)|_p^{-[L:\qp]}.\]
\end{prop}

\begin{proof}  First observe that using Lemma \ref{tensor}, property h) and the fact that $\O$ is regular we obtain canonical isomorphisms 
\begin{eqnarray*}
\xymatrix{
  {\u_\O} \ar[rr]^<(0.3){\O^n\otimes_{\La(G)}a} & &{\d_\O(\O^n\otimes^\mathbb{L}_{\La(G)}\cA)}   }&\cong & \d_\O(\O\otimes^\mathbb{L}_{\La_\O(G)}\cA(\rho^*)))\\
&\cong& \d_\O(\H^\cdot_\B(G,\cA(\rho^*)))\\
&\cong& \prod_{i\in\mathbb{Z}}\d_\O\big(\H^i_\B(G,\cA(\rho^*))\big)^{(-1)^i}.
\end{eqnarray*}
After tensoring with $L$ and identifying then the factors at the end with the unit object by acyclicity we recover the definition of the leading term $(\cA,a)^*(\rho).$  From this description it becomes clear that the latter can also be written modulo $\O^\times$ as  the product over $i$ of the following maps \[ \xymatrix{
  (\u_\O)_L \ar[r]^<(0.2){def} & {\d_\O(\H^i_\B(\cA(\rho^*))\big)_L^{(-1)^i}} \ar[r]^<(0.3){acyc} & {\u_L}   }.\] The result follows now from Remark \ref{eltdivisor}.
\end{proof}

\begin{remark}{\em In the special case that the complex $\cA(\rho^*)\otimes_{\bz_p}\bq_p$ is acyclic, the leading term  $\L^*(\rho)$ is just the value of $\L$ at $\rho$ (in the sense of Remark \ref{value}) and so the result of
Proposition \ref{gecp} recovers the results of \cite[thm.\
3.6]{cfksv}, \cite[prop.\ 6.3 ]{ven-crelle}  and \cite[rem.\
4.1.13]{fukaya-kato}.}\end{remark}

\section{Global Zeta isomorphisms}

In this section we recall the Tamagawa Number Conjecture as formulated by Fukaya and Kato in \cite{fukaya-kato}. 
%In addition, at the end of the section we indicate how the approach of
%\S\ref{tgc}   leads to an interesting refinement of the
%observations made by Fukaya and Kato in \cite[\S2.5]{fukaya-kato}.

\subsection{Galois cohomology}\label{appendix}

The main reference for this section is \cite[\S
1.6]{fukaya-kato}, but see also \cite{bufl01}, here we use the same notation as in \cite{ven-BSD}. For simplicity we
assume that $p$ is odd throughout this section. Let
$U={\rm spec}(\z[\frac{1}{S}])$ be a dense open subset of ${\rm spec}(\z)$
where $S$ contains $S_p:=\{p\}$ and $S_\infty:=\{\infty\}$ (by abuse of
notation). We write $G_S$ for the Galois group of the maximal
outside $S$ unramified extension of $\Q.$ Let $X$ be a topological
abelian group with a continuous action of $G_S.$ 
As usual we write $\r(U,X)$ ($\r_c(U,X)$) for global Galois cohomology with restricted ramification (and compact support), for any place $v$ of $\Q$ we denote by $\r(\Q_v,X)$ local Galois cohomology. Let $L$ be a finite extension of $\qp$ with ring of integers $\O$ and  $V$ a
finite dimensional $L$-vector space  with continuous
$G_{\Q_v}$- and $G_\Q$-action, respectively. Then the {\em finite parts} of global and local cohomology are  written as $\r_f(\Q,V)$ and $\r_f(\Q_v,V),$ respectively. There is a canonical exact triangle
 \be \label{trian-c-f} \xymatrix{
  {\r_c(U,V)} \ar[r]^{ } & {\r_f(\Q,V)} \ar[r]^{ } & {\bigoplus_{S}\mathrm{R\Gamma}_f(\Q_v,V) }
  \ar[r]^{ } &  .   }\ee

%\subsubsection{Duality}

%Let $G,$ $\La=\La(G),$ $\T$ as in section \ref{ETNC}. By abuse of notation we write $-^*$ for both (derived) functors $\rhom_\La(-,\La)$ and $\rhom_{\La^\circ}(-,\La^\circ).$ Then Artin-Verdier/Poitou-Tate duality induces the existence of the following distinguished triangle in the derived category of $\La$-modules
%\be \label{AVPT} \xymatrix@C=0.5cm{
%   { \r_c(U,\T) }\ar[rr]^{ } && {\r(U,\T^*(1))^*[-3] }\ar[rr]^{ } && {\T^+} \ar[r] &   }\ee
%   and similarly for $T$ (a Galois stable $\O$-lattice of $V$) and $V$ as coefficients (with   $\La$ replaced by $\O $ and  $L,$ respectively).\footnote{A more precise form to state the duality is the following. Let $\r_{(c)}(U,\T)$ be defined like $\r_{c}(U,\T)$  but using Tate cohomology $\widehat{\r}  (\R,\T)$ instead of the usual group cohomology $\r(\R,\T).$ Then one has isomorphisms \[ \r(U,\T^*(1))^*\cong\r_{(c)}(U,\T)[3]\cong \r(U,\T^\vee(1))^\vee\] where $-^\vee=\Hom_{cont}(-\qp/\zp)$ denotes the Pontryagin dual.}

%For the finite parts one obtains from Artin-Verdier/Poitou-Tate  and local Tate-duality the following isomorphisms
%\be \label{local-finite-dual}\r_f(\ql,V)\cong (\r(\ql,V^*(1))/\r_f(\ql,V^*(1)))^*[-2], \ee
%\be \label{global-finite-dual}\r_f(\Q,V)\cong \r_f(\Q,V^*(1))^*[-3].\ee

For any prime $\ell$ set $t_\ell(V):=D_{dR}(V)/D_{dR}^0(V)$ if $\ell=p$ and $t_\ell(V)=0$ otherwise. For the definition of the canonical isomorphism in $\C_L$
\begin{eqnarray}\label{eta}
\eta_\ell(V):& \u_{L} &\to
 \d_{L}(\mathrm{R\Gamma}_f(\Q_l,V))\d_L(t_\ell(V)).
\end{eqnarray}
we refer the reader to \cite[\S 2.4.4]{fukaya-kato} or the appendix of \cite{ven-BSD}.

\subsection{$K$-Motives over $\Q$}

 %We shall simply view motives in the naive sense, as being defined by a
%collection of realizations satisfying certain axioms, together
%with their motivic cohomology groups. The archetypical motive is
%$h^i(X)$ for a smooth projective variety $X$ over $\Q$ with its
%obvious {\'e}tale cohomology
%$\H^i_{\acute{e}t}(X\times_\Q\overline{\Q},\Q_l)$, singular
%cohomology $\H^i(X(\mathbb{C}),\Q)$ and de Rham cohomology
%$\H^i_{dR}(X/\Q),$ their additional structures and comparison
%isomorphisms. More general, let $K$ be a finite extension of $\Q.$
%A $K$-motive $M$ over $\Q,$ i.e.\ a motive over $\Q$ with an
%action of $K,$ will be given by the following data, which for
%$M=h^n(X)_K$ arise by tensoring the above cohomology groups by $K$
%over $\Q :$

For the background on this material we refer the reader to \cite[\S 2.2, 2.4]{fukaya-kato}, \cite[\S 3]{bufl01} or the survey \cite[\S 2]{ven-BSD}. Let $K$ be a finite extension of $\Q.$ As usual we write $M_B, M_{dR}$ and $M_{\lambda}$ for the Betti, de Rham and $\lambda$-adic realizations of a $K$-motive $M$ where $\lambda$ varies over all finite places of $K.$ Also, we write $t_M=M_{dR}/M^0_{dR}$ for the tangent space of $M.$ Henceforth, for any commutative ring $R$
  and $R[G(\bbC/\R)]$-module $X$ we denote by $X^+$ and $X^-$ the
$R$-submodule of $X$ on which complex conjugation  acts by $+1$ and $-1,$
respectively.
In our later calculations we will use the following isomorphisms

\begin{itemize}
 \item[$\bullet$] The comparison between the Betti- and the $\lambda$-adic realization induces  canonical isomorphisms of $K_\lambda$ and $K_l$-modules
\begin{equation}\label{B-l}
g_\lambda^+:K_\lambda\otimes_K M_B^+ \cong M_\lambda^+     \mbox{
and } g_\ell^+:K_\ell\otimes_K M_B^+ \cong M_\ell^+.
\end{equation}
\item[$\bullet$] The comparison between de Rahm and Betti-cohomology induces the
 the $\R$-linear period map
\begin{equation}\label{alphaM}
 \xymatrix{
{\R\otimes_\Q M_B^+} \ar[r]^{\alpha_M} & {\R\otimes_\Q
t_M.} }
\end{equation}
%We say that $M$ is {\em critical}  if this map is an 
%isomorphism.
 \item[$\bullet$]
The comparison between $p$-adic and de Rham cohomology induces an
isomorphism of
$K_\lambda$-vector spaces
\begin{equation}
 \xymatrix{
{  t_p(M_\lambda)=D_{dR}(M_\lambda)/D_{dR}^0(M_\lambda)}
\ar[r]^<(0.2){g_{dR}^t}_<(0.2){\cong} & { K_\lambda\otimes_K t_M.}
}
\end{equation}
\end{itemize}

The {\em motivic cohomology} $K$-vector spaces  $\H_f^0(M):=\H^0(M)$ and
$\H^1_f(M)$ may be defined by algebraic $K$-theory or motivic
cohomology a la Voevodsky. They are conjectured to be finite
dimensional. For the    definition we refer the reader to \cite{bufl01}.

%\begin{ex}\label{ex-mot-coh}
%A) For the Tate motive we have $\H_f^0(\Q(1))=\H_f^1(\Q(1))=0$ and
%for its Kummer dual $\H_f^0(\Q)=\Q$ while $\H_f^1(\Q)=0.$
%
%B) If $M=h^1(A)(1)$  for an abelian variety $A$ over $\Q$ one has $\H^0_f(M)=0$ and $\H^1_f(M)=A^\vee(\Q)\otimes_\z\Q.$
%
%C) For $M=h^0(spec(F))$  we have $\H_f^0(M)=\Q$ and $\H_f^1(M)=0$
%while for $M^*(1)=h^0(spec(F))(1)$ one has $\H_f^0(M^*(1))=0$ and
%$\H_f^1(M^*(1))= \O_F^\times\otimes_\z\Q.$ More general, for an
%$K$-Artin motive $[\rho]$ one has $\H^0_f([\rho])=K^n,$ where $n$
%is the multiplicity with which $\Q$ occurs in $[\rho].$
%
%Unfortunately the functor $\H^i_f$ does not behave well with
%tensor  products, i.e.\ in general one cannot derive
%$\H^*_f([\rho]\otimes h^1(A)(1))$ from $\H^*_1([\rho])$ and
%$\H^*_1(h^1(A)(1))$ (e.g.\ in form of a K\"unneth formula).
%\end{ex}

\subsection{The Tamagawa Number Conjecture}

%In \cite{BK} Bloch and Kato formulated a vast generalization of
%the analytic class number formula and the BSD-conjecture. While
%the conjecture of Deligne and Beilinson links the order of
%vanishing of the $L$-function attached to a motive $M$ to its
%motivic cohomology and claims rationality of special $L$-values or
%more general  leading coefficients (up to periods and regulators)
%the Tamagawa number conjecture by Bloch and Kato predicts the
%precise $L$-value in terms of Galois cohomology (assuming the
%conjecture of Deligne-Beilinson).
%
%Later, Fontaine and Perrin-Riou \cite{fp} found an equivalent
%formulation using (commutative) determinants instead of (Tamagawa)
%measures\footnote{The name comes from an analogy with the theory
%of algebraic groups, see \cite{BK}.}. In this section we follow
%closely their approach.

Let us first recall from \cite[\S 2.2.2]{fukaya-kato} that, for each embedding $K\to\mathbb{C},$ the complex $L$-function
attached to a $K$-motive $M$ is defined as Euler product \[L_K(M,s)=\prod_p P_p(M,p^{-s})^{-1}\]
  for the real part of $s$ large enough. We assume meromorphic continuation and    write as usual $L^*_K(M)\in \mathbb{R}^\times$ and $r(M) \in\mathbb{Z}$ for its leading coefficient and order of vanishing at $s=0,$ respectively.

To establish a link between $L^*(M)$ and Galois cohomology one uses the   {\em fundamental line}:

\begin{eqnarray*}
\Delta_K(M):&=&\d_K(\H^0_f(M))^{-1}\d_K(\H^1_f(M))\d_K(\H^0_f(M^*(1))^*)
\d_K(\H^1_f(M^*(1))^*)^{-1}\\ && \d_K(M_B^+)\d_K(t_M)^{-1}.
\end{eqnarray*}

The conjecture formulated in  \cite[\S 2.2.7]{fukaya-kato} induces a canonical isomorphism in $\C_{K_\br}$
(period-regulator map)
\begin{equation}
\vartheta_\infty: K_\bbR\otimes_K\Delta_K(M)\cong \u_{K_\bbR}.
\end{equation}

In addition, a standard conjecture  on Cycle class and Chern class maps induces, for any place $\lambda$ above $p,$ the $p$-adic period-regulator isomorphism in $C_{K_\lambda}$ (involving the maps $\eta_l$)
%($p$-adic period-regulator map)
%\begin{equation}\label{vartheta}
%\vartheta_p:  \Delta_K(M)_{K_p}\cong
%\d_{K_p}({\mathrm{R\Gamma}_c(U,M_p)})^{-1} ,
%\end{equation}
% which induces for any place $\lambda$ above $p$
\begin{equation}
\vartheta_\lambda:  \Delta_K(M)_{K_\lambda}\cong
\d_{K_\lambda}({\mathrm{R\Gamma}_c(U,M_\lambda)})^{-1} .
\end{equation}

Now let $F_\infty$ be a $p$-adic Lie extension of $\Q$ with Galois group
$G=G(F_\infty/\Q).$ By $\La=\La(G)$ we  denote its Iwasawa algebra. For a
$\Q$-motive $M$ over $\Q$ we fix a $G_\Q$-stable $\zp$-lattice $T_p$ of
$M_p$ and define a left $\La$-module
\[\T :=\La\otimes_\zp T_p\]
on which  $\La$ acts via multiplication on the left factor from
the left  while $G_\Q$ acts diagonally via $g(x\otimes y)=x
\bar{g}^{-1}\otimes g(y),$ where $\bar{g}$ denotes the image of
$g\in G_\Q$ in $G.$  

Let  $\lambda$ a finite place of
$K,$ $\O_\lambda$ the ring of integers of the completion
$K_\lambda$ of $K$ at $\lambda$ and assume that $\rho:G\to
GL_n(\O_{\lambda})$ is a continuous representation of $G$ which,
  for some suitable choice of a basis, is the $\lambda$-adic
realization $N_\lambda$ of a  $K$-motive $N.$ We also write
$\rho$ for the induced ring homomorphism $\La(G)\to M_n(
\O_{\lambda})$ and we consider $\O_\lambda^n$ as a right
$\Lambda(G)$-module via action by $\rho^t$ on the left,  viewing
$\O_\lambda^n$ as set of column vectors (contained in
$K_\lambda^n.)$ Note that, setting $M(\rho^*):= N^*\otimes M,$ we
obtain an isomorphism of Galois representations
\[\O_\lambda^n\otimes_{\La(G)} \T\cong T_\lambda(M(\rho^*)),
\]
where $T_\lambda(\rho^*)$ is the $\O_\lambda$-lattice
$\rho^*\otimes T_p$ of $M(\rho^*)_\lambda$ and $\rho^*$ denotes
the contragredient(=dual) representation of $\rho.$

\begin{conj}[{Fukaya/Kato \cite[conj.\ 2.3.2]{fukaya-kato} or \cite[conj.\ 4.1 ]{ven-BSD}}]\label{equivintegrality}
There exists a (unique)
 isomorphism in $\C_\La$ \[\zeta_\La(M):=\zeta_\La(\T): \u_\La \to \d_\La(\r_c(U,\T))^{-1}\]
with the following property: for all $K, \lambda$ and $ \rho$ as above the  (generalized) base
change $K_\lambda^n\otimes_\La -$ sends $\zeta_\La(M)$ to
%\[\zeta_{\O_\lambda}(T_\lambda(\rho^*)):\u_{\O_\lambda} \to \d_{\O_\lambda}(\mathrm{R\Gamma}_c(U,T_\lambda(\rho^*)))^{-1},\]
%
%which in turn induces via $K_\lambda\otimes_{\O_\lambda}-$ the  
% map ${\zeta_{\O_\lambda}(T_\lambda(\rho^*))_{K_\lambda}:}$
 \[\label{zetaintegral}\xymatrix{
   {\u_{K_\lambda}} \ar[rr]^<(0.3){\zeta_K(M(\rho^*))_{K_\lambda}} &
    & {\Delta_K(M(\rho^*))_{K_\lambda}} \ar[r]^<(0.15){\vartheta_\lambda} &   {\d_{K_\lambda}({\mathrm{R\Gamma}_c(U,M(\rho^*)_\lambda)})^{-1},}
       }\]
       where  $\zeta_K(M(\rho^*)):\u_K\to\Delta_K(M(\rho^*))$ is the unique isomorphism such that, for every embedding   $K \to \bbC,$ the leading coefficient $L_K^*(M(\rho^*))$  is equal to
the composite
 \[ \xymatrix{
   {\phantom{L_K^*(M(\rho^*)}\u_\bbC}\ar[rr]^{\zeta_K(M(\rho^*))_\bbC} &   & {\Delta_K(M(\rho^*))_\bbC} \ar[rr]^{(\vartheta_\infty)_\bbC} &  & {\u_\bbC}.   }\]
\end{conj}

%Note that this conjecture assumes conjecture \ref{integrality} for
%all $K$-motives $M(\rho^*)$ with varying $K.$ Furthermore, it is
%independent of the choice of $S$ and of the lattices $T_p(M)$ and
%$T_\lambda(\rho^*).$

It is easily shown that
this conjecture implies the Equivariant Tamagawa Number Conjecture formulated
by Flach and the first named author in \cite[conj.\ 4(iv)]{bufl01} and hence also implies the `main conjecture of
non-abelian Iwasawa theory' discussed by Huber and Kings in \cite{hu-ki}.

\section{The interpolation formula for Tate motives}\label{Q(1)}

In this section we shall give a first explicit application of the
formalism developed in \S\ref{bock}. More precisely, we show that the
`$p$-adic Stark conjecture at $s=1$', as formulated by Serre in
\cite{ser} and discussed by Tate in \cite[chap. VI, \S5]{tate},
can be naturally interpreted as an interpolation formula for the
leading term (in the sense of Definition \ref{leadingterm-rho}) of
 certain global Zeta isomorphisms that are predicted to exist by Conjecture \ref{equivintegrality} in terms
of the leading terms (in the classical sense) of suitable $p$-adic
Artin $L$-functions.  Interested
readers can find further explicit results concerning Conjecture \ref{equivintegrality} in the
special case that we consider here in both   \cite{br-bu2} and
\cite{bufl05}.

Throughout this section we fix an odd prime $p$ and a totally real
Galois extension $F_\infty$ of $\bq$ which contains the cyclotomic
$\bz_p$-extension $\bq_{\rm cyc}$ of $\bq$ and is such that $G :=
\mathrm{G}(F_\infty/\bq)$ is a compact $p$-adic Lie group. We assume further
that $F_\infty/\bq$ is unramified outside a finite set of prime numbers
 $S$ (which therefore contains $p$). We
set $H := \mathrm{G}(F_\infty/\bq_{\rm cyc})$ and $\Gamma := \mathrm{G}(\bq_{\rm cyc}/\bq)
\cong G/H$. We fix a subfield $E$ of $F_\infty$ which has finite degree
over $\bq$ and, for simplicity, we assume throughout that the
following condition is satisfied:
\be\label{disjoint} E \cap \bq_{\rm cyc} = \bq \mbox{ and }  E_w\cap \bq_{p, \rm cyc} =\bq_p \mbox{  for
all }w\in S_p(E).\ee

We let $\T$ denote the $\Lambda(G)$-module $\Lambda(G)$ endowed
with the following action of $\mathrm{G}(\overline{\bq}/\bq)$: each $\sigma \in
\mathrm{G}(\overline{\bq}/\bq)$ acts on $\T$ as multiplication by the element
$\chi_{\rm cyc}(\bar{\sigma})\bar{\sigma}^{-1}$ where
$\bar{\sigma}$ denotes the image of $\sigma$ in $G$ and $\chi_{\rm
cyc}$ is the cyclotomic character $G \to \Gamma \to \bz_p^\times$.
For each subfield $F$ of $F_\infty$ which is Galois over $\bq$ we let
$\T_F$ denote the $\La(\mathrm{G}(F/\bq))$-module $
\La(\mathrm{G}(F/\bq))\otimes_{\La (G)}\T$. We also set $U := \Spec (\bz
 [\frac{1}{S}])$ and note that for each such field $F$ there is a
natural isomorphism in $D^p(\La (\mathrm{G}(F/\bq)))$ of the form
\be\label{descent} \La(\mathrm{G}(F/\bq))\otimes^{\mathbb L}_{\La
(G)}R\Gamma_c(U,\T) \cong R\Gamma_c(U,\T_F).\ee

We regard each character of $\bar{G}:=G(E/\bq)$ as a character of $G$ via the
natural projection $G \to \bar{G}$. For any field $C$ we write
$R_C(\bar{G})$ for the ring of finite dimensional $C$-valued
 characters of $\bar{G}$, and for each $C$-valued character $\rho$ we fix a
representation space $V_\rho$ of character $\rho$. For any $\bq_p
[\bar{G}]$-module $N$, resp. endomorphism $\alpha$ of a $\bq_p
[\bar{G}]$-module $N$, we
 write $N^\rho$ for the $\bc_p$-module
 \[\Hom_{\bar{G}}(V_\rho,\bc_p\otimes_{\bq_p}N)\cong ((V_{\rho^*})_{\bc_p}\otimes_{\bq_p} N)_{\bar{G}},\] resp.\ $\alpha^\rho$ for the induced endomorphism of
 $N^\rho$. We use similar notation for complex characters $\rho$ and $\bq [\bar{G}]$-modules $N$.

For any abelian group $A$ we write $A\hat\otimes \bz_p$ for its
$p$-adic completion $\varprojlim_n A/p^nA$.

\subsection{Leopoldt's Conjecture.} We recall that Leopoldt's
conjecture (for the field $E$ at the prime $p$) is equivalent to
the injectivity of the natural localisation map
\[ \lambda_p: \O_E\left[{1\over {p}}\right]^\times\otimes\bz_p \to \prod_{w \in
S_p(E)}E_w^\times\hat\otimes\bz_p,\]
where $S_p(E)$ denotes the set of places of $E$ which lie above
$p$. If $\rho \in R_{\bc_p}(\bar{G})$, then we say that Leopoldt's
conjecture `is valid at $\rho$' if
$(\bq_p\otimes_{\bz_p}\ker(\lambda_p))^\rho = 0$.

We set $c_\gamma := c_{\chi_{\rm cyc},\gamma} \in \bq_p^\times$ and for
 each $\rho \in R_{\bc_p}(\bar{G})$ we define
\[ \langle \rho, 1\rangle := \dim_{\bc_p}(H^0(\bar{G},V_\rho)) = \dim_{\bc_p}((\bq_p)^\rho).\]
\begin{lem}\label{exp-des} We fix $\rho \in R_{\bc_p}(\bar{G})$ and assume that Leopoldt's conjecture is valid at $\rho$.

\begin{itemize}
\item[(i)] There are canonical isomorphisms
\[ (H^i_c(U,\T_E)\otimes_{\bz_p}\bq_p)^\rho \cong \begin{cases}
(\cok(\lambda_p)\otimes_{\bz_p}\bq_p)^\rho, &\text{if $i = 2$}\\
(\bq_p)^\rho,&\text{if $i = 3$}\\
0, &\text{otherwise.}\end{cases}\]
\item[(ii)] $R\Gamma_c(U,\T)\in \Sigma_{{\rm ss}-\rho}$ and
$r_G(R\Gamma_c(U,\T))(\rho) = \langle \rho, 1\rangle$.
\item[(iii)] For each $w \in S_p(E)$ we write ${\rm
N}_{E_w/\bq_p}$ for the morphism $E_w^\times\hat\otimes\bz_p \to
\bq_p^\times\hat\otimes \bz_p$ that is induced by the field
theoretic norm map. Then, with respect to the identifications
given in claim (i), the morphism
\[ (\B^2)^\rho: (H^2_c(U,\T_E)\otimes_{\bz_p}\bq_p)^\rho \to (H^3_c(U,\T_E)\otimes_{\bz_p}\bq_p)^\rho, \]
is induced by the map
\[\log_{p,\gamma,E}: \prod_{w \in S_p(E)}E_w^\times\hat\otimes \bz_p \to \bz_p\]
which sends each element $(e_w)_w$ to $-\sum_w
c_\gamma^{-1}\log_p({\rm N}_{E_w/\bq_p}(e_w))$.
\end{itemize}
\end{lem}
\begin{proof} Claim (i) can be verified by combining the exact
cohomology sequence of the tautological exact triangle
\be\label{dist-tri} R\Gamma_c(U,\T_E) \to R\Gamma(U,\T_E) \to
\bigoplus_{\ell \in S}R\Gamma(\bq_\ell,\T_E)\to
R\Gamma_c(U,\T_E)[1]\ee
together with the canonical identifications $H^i(U,\T_E) \cong
H^i(\O_{E}[\frac{1}{S}],\bz_p(1))$ and $H^i(\bq_\ell,\T_E) \cong
\bigoplus_{w \in S_\ell(E)}H^i(E_w,\bz_p(1))$ and an explicit
computation of each of the groups
$H^i(\O_{E}[\frac{1}{S}],\bz_p(1))$ and $H^i(E_w,\bz_p(1))$.
 As this is routine we leave explicit details to the reader
except to note that $H^2_c(U,\T_E)\otimes_{\bz_p}\bq_p$ is
canonically isomorphic to $\cok(\lambda_p)\otimes_{\bz_p}\bq_p$
(independently of Leopoldt's conjecture), whilst the fact that $E$
is totally real implies the vanishing of
$(H^1_c(U,\T_E)\otimes_{\bz_p}\bq_p)^\rho$ is equivalent to that
of $(\ker(\lambda_p)\otimes_{\bz_p}\bq_p)^\rho$.

To prove claim (ii) and (iii) we write $E_{\rm cyc}$, $E_{w,\rm cyc}$
for each $w \in S_p(E)$ and $\bq_{p,\rm cyc}$ for the cyclotomic
$\bz_p$-extensions of $E$, $E_w$ and $\bq_p$. Then
(\ref{disjoint}) induces a direct product decomposition
$\mathrm{G}(E_{\rm cyc}/\bq) \cong \Gamma\times \bar{G}$ and hence allows us to
identify $\Gamma$ with each of $\mathrm{G}(E_{\rm cyc}/E)$,
$\mathrm{G}(E_{w,\rm cyc}/E_w)$ and $\mathrm{G}(\bq_{p,\rm cyc}/\bq_p)$. We note
also that, in terms of the notation of Lemma \ref{tensor}, the
isomorphism
 (\ref{descent}) (with $F = E_{\rm cyc}$) induces a canonical isomorphism in $D^p(\La (\Gamma))$ of the form
\be\label{reinter} R\Gamma_c(U,\T)_\rho \cong
\O^n\otimes_{\bz_p[\bar{G}]}R\Gamma_c(U,\T_{E_{\rm cyc}}),\ee
where $\Gamma$ acts naturally on the right hand factor in the
tensor product.

From Lemma \ref{tensor}(i), we may therefore deduce that
$R\Gamma_c(U,\T) \in \Sigma_{{\rm ss}-\rho}$ if and only if
$\O^n\otimes_{\bz_p[\bar{G}]}R\Gamma_c(U,\T_{E_{\rm cyc}})\in
\Sigma_{\rm ss}$. But the latter containment can be easily
verified by using the criterion of Remark \ref{ss-remark}(ii):
indeed, one need only note that $H^i_c(U,\T_{E_{\rm cyc}})$ is finite
if $i \notin \{2,3\}$, that $H^3_c(U,\T_{E_{\rm cyc}})$ identifies
with $\bz_p$ (as a $\Gamma$-module) and that the exact sequences
of (\ref{bock-decomp}) combine with the descriptions of claim (i)
to imply that
$((\bq_p\otimes_{\bz_p}H^1_c(U,\T_{E_{\rm cyc}}))^\rho)^\Gamma$ and
$(\bq_p\otimes_{\bz_p}H^1_c(U,\T_{E_{\rm cyc}}))^\rho_{ \Gamma}$ both
vanish. In addition, the same observations combine with Lemma
\ref{tensor}(ii) to imply that $r_G(R\Gamma_c(U,\T))(\rho) =
\dim_{\bc_p}((\bq_p)^\rho)$.

Regarding claim (iii), the isomorphism (\ref{reinter}) combines
with Lemma \ref{tensor}(i) and Lemma \ref{bockstein} to imply that
$(\B^2)^\rho = \hat (\B^2)^\rho$ where $\hat \B^2$ is the Bockstein
morphism of the pair $(R\Gamma_c(U,\T_{E_{\rm cyc}}), T)$. Also, the
cohomology sequence of the triangle (\ref{dist-tri}) gives rise to
a commutative diagram
\[ \begin{CD} \bigoplus_{w \in S_p(E)}\bq_p\otimes_{\bz_p}H^1(E_w,\bz_p(1)) @> >>
\bq_p\otimes_{\bz_p}H^2_c(U,\T_E)\\
@V (-1) \times (\bq_p\otimes_{\bz_p}\B^1_w)_wVV @VV \bq_p\otimes_{\bz_p}\hat \B^2V\\
\bigoplus_{w \in S_p(E)}\bq_p\otimes_{\bz_p}H^2(E_w,\bz_p(1)) @>
>> \bq_p\otimes_{\bz_p}H^3_c(U,\T_E)\end{CD}\]
where the upper row is the (tautological) surjection induced by
the canonical identifications $H^1(E_w,\bz_p(1)) \cong
E_w^\times\hat\otimes \bz_p$ and
$\bq_p\otimes_{\bz_p}H^2_c(U,\T_E) \cong
\bq_p\otimes_{\bz_p}\cok(\lambda_p)$, the lower row is the
surjection induced by the canonical identifications
$H^2(E_w,\bz_p(1)) \cong \bz_p$ and $H^3_c(U,\T_E) \cong \bz_p$
together with the identity map on $\bz_p$, $\B^1_w$ is the
Bockstein morphism of the pair $(R\Gamma(E_{w,\rm cyc},\bz_p(1)), T)$
and the factor $-1$ occurs on the left hand vertical arrow because
 the diagram compares Bockstein morphisms in degrees $1$ and $2$.

Further, for each $w \in S_p(E)$ the natural isomorphism (in
$D^p(\bz_p)$)
\[ \bz_p\otimes^{\mathbb L}_{\bz_p[\mathrm{G}(E_w/\bq_p)]}R\Gamma(E_w,
\bz_p(1)) \cong R\Gamma(\bq_p,\bz_p(1))\]
induces a commutative diagram
\[ \begin{CD} H^1(E_w,\bz_p(1)) @> >>
H^1(\bq_p,\bz_p(1))\\
@V \B^1_w VV @VV \B^1_pV\\
H^2(E_w,\bz_p(1)) @> >> H^2(\bq_p,\bz_p(1))
\end{CD}\]
where the upper horizontal arrow is induced by the canonical
identifications $H^1(E_w,\bz_p(1)) \cong
E_w^\times\hat\otimes\bz_p$ and $H^1(\bq_p,\bz_p(1)) \cong
\bq_p^\times\hat\otimes\bz_p$ together with the map ${\rm
N}_{E_w/\bq_p}$, the lower horizontal arrow is induced by the
canonical identifications $H^2(E_w,\bz_p(1)) \cong \bz_p$ and
 $H^2(\bq_p,\bz_p(1)) \cong \bz_p$ together with the identity map on $\bz_p$, and $\B^1_p$ is
the Bockstein morphism of the pair
$(R\Gamma(\bq_{p,\rm cyc},\bz_p(1))),T)$.
To prove claim (iii) it thus suffices to recall that, with respect
to the identifications $H^1(\bq_p,\bz_p(1)) \cong
\bq_p^\times\hat\otimes\bz_p$ and $H^2(\bq_p,\bz_p(1)) \cong
\bz_p$, the map $\B^1_p$ is equal to $c_\gamma^{-1}\cdot\log_p$
(see, for example, \cite[p.\ 352]{burns-greither2003}).
\end{proof}

\subsection{The $p$-adic Stark conjecture at $s =
1$.}\label{pasc1} For each character $\chi \in R_\bc(\bar{G})$ we
write $L_S(s,\chi)$ for the Artin $L$-function of $\chi$ that is
truncated by removing the Euler factors attached to primes in $S$
(cf. \cite[Chap. 0, \S4]{tate}). Then, for each character $\rho
\in R_{\bc_p}(G)$ there exists a unique $p$-adic meromorphic
function $L_{p,S}(\cdot, \rho): \bz_p \to \bc_p$ such that for
each strictly negative integer $n$ and each isomorphism $j:\bc
_p\cong \bc$ one has
\[ L_{p,S}(n,\rho)^j = L_S(n,(\rho\cdot\omega^{n-1})^j )\]
where $\omega: \mathrm{G}(\overline{\bq}/\bq) \to \bz_p^\times$ is the
Teichm\"uller character. Henceforth we will fix such an isomorphism $j$ and often omit it from the notation.  This function is the `($S$-truncated)
$p$-adic Artin $L$-function' of $\rho$ that is constructed by
Greenberg in \cite{greenberg} by combining techniques of Brauer
induction with the fundamental results of Deligne and Ribet
\cite{deligne-ribet} and Cassou-Nogu\`es \cite{c-n}.

In this section we recall a conjecture of Serre regarding the
`leading term at $s=1$' of $ L_{p,S}(s,\rho)$. To this end we set
$E_\infty := \br \otimes_\bq E \cong \prod_{\Hom(E,\bc)}\br$ and
write $\log_\infty(\O_E^\times)$ for the inverse image of
$\O_E^\times \hookrightarrow E_\infty^\times$ under the
 (componentwise) exponential map $\exp_\infty: E_\infty \to E_\infty^\times.$ We set $E_0 := \{ x \in E: {\rm
Tr}_{E/\bq}(x) = 0\}$. Then $\log_\infty(\O_E^\times)$ is a
lattice in $\br \otimes_\bq E_0$ and so there is a canonical
 isomorphism of $\bc [\bar{G}]$-modules $\mu_\infty: \bc \otimes_\bz\log_\infty(\O_E^\times) \cong \bc\otimes_\bq E_0$.
 This implies that the $\bq[\bar{G}]$-modules $E_0$ and $\O_{E}^\times\otimes _\bz \bq$ are (non-canonically) isomorphic.
 We note also that the composite
morphism
\begin{multline}\label{inf-comp} \log_\infty(\O_E^\times) \xrightarrow{\rm exp_\infty}
\O_E^\times \xrightarrow{\lambda_p} \prod_{w \in
S_p(E)}U^1_{E_w}\\
\xrightarrow{(u_w)_w \mapsto (\log_p(u_w))_w} \prod_{w \in
S_p(E)}E_w \cong \bq_p\otimes_{\bq}E,\end{multline}
factors through the inclusion $\bq_p\otimes_{\bq}E_0 \subset
\bq_p\otimes_{\bq}E$ and hence induces an isomorphism of
$\bq_p[\bar{G}]$-modules $\mu_p: \bq_p
\otimes_\bz\log_\infty(\O_E^\times) \cong \bq_p\otimes_\bq E_0$.

\begin{conj}\label{pasc} For each $\rho \in R_{\bc_p}(\bar{G})$ we set
\[ L_{p,S}^*(1,\rho) := \lim_{s \to 1}(s-1)^{\langle \rho, 1\rangle}\cdot L_{p,S}(s,\rho).\]
Then
$L^*_{p,S}(1,\rho)$ is  equal to the leading term of $L_{p,S}(s,\rho)$ at $s =
1$, and   for   each choice of isomorphism of $\bq [\bar{G}]$-modules $g: E_0 \to
\O_{E}^\times\otimes _\bz \bq$  one has
\[ \frac{L_{p,S}^*(1,\rho) }{{\rm det}_{\bc_p}((\bc_p\otimes_{\bq_p}\mu_p)\circ (\bc_p
\otimes_\bq g))^\rho} = \frac{L_S^*(1,\rho )}{{\rm
det}_\bc (\mu_\infty\circ (\bc \otimes_\bq
g))^{\rho }}.\]
\end{conj}

\begin{remark} {\em This conjecture is the `$p$-adic Stark conjecture at $s =1$' as discussed by Tate in
\cite[Chap. VI, \S5]{tate} (where it is attributed to
 Serre \cite{ser}). More precisely, there are some imprecisions in the statement of \cite[Chap. VI, \S5]{tate} (for example,
 and as already noted by Solomon
in \cite[\S3.3]{solo}, the intended meaning of the symbols `$\log
U$' and `$\mu_p$' in \cite[p.\ 137]{tate} is unclear) and
Conjecture \ref{pasc} represents a natural clarification of the
presentation given in loc.\ cit.\ }
\end{remark}

\begin{remark}\label{colmez} {\em We fix a subgroup $J$ of $\bar{G}$ and
 write ${\rm 1}_J$ for the trivial character of $J$. If $\rho = \Ind_J^{\bar{G}}{\rm 1}_J$, then the inductive behaviour of $L$-functions
 combines with the analytic class number formula for $E^J$ to show that
 Conjecture \ref{pasc} is valid for $\rho$ if and only if the residue at $s=1$ of the $p$-adic zeta function of
 the field $E^J$ is non-zero and equal to $2^{[E^J:\bq]-1}hR_pe_p/\sqrt{|d|}$
  where $h, R_p$ and $d$ are the class number, $p$-adic regulator and absolute discriminant of
  $E^J$ respectively and $e_p := \prod_{v \in S_p(E^J)}(1 - {\rm N}v^{-1})$ (cf.\ \cite[rem., p.\ 138]{tate}). From the
  main result of Colmez in \cite{col} one  may thus deduce that
   Conjecture \ref{pasc} is valid for $\rho = \Ind_J^G{\rm 1}_J$ if and only if Leopoldt's conjecture is valid for $E^J$.
   We note also that if Leopoldt's conjecture is valid for $E$,
then it is valid for all such intermediate fields $E^J$.}
\end{remark}

\subsection{The interpolation formula} We now reinterpret the equality of Conjecture \ref{pasc} as an interpolation formula for the Zeta isomorphism $\zeta_{\Lambda(G)}(\T).$ 

\begin{thm}\label{st-thm} If Conjecture \ref{pasc} is valid, then for each $\rho\in R_{\bc_p}(\bar{G})$ one has
$R\Gamma_c(U,\T) \in \Sigma_{{\rm ss}-\rho}$,
$r_G(R\Gamma_c(U,\T))(\rho) = \langle \rho, 1\rangle$ and
\be\label{nice} c_\gamma^{\langle \rho, 1\rangle} \cdot
\zeta_{\Lambda(G)}(\T)^*(\rho) = L_{p,S}^*(1,\rho).\ee
\end{thm}

\begin{remark}\label{equa-lead}{\em One can naturally interpret (\ref{nice}) as an equality of leading terms of $p$-adic meromorphic functions.
 Indeed, whilst Conjecture \ref{pasc} predicts  that $L_{p,S}^*(1,\rho)$ is the leading term at $ s= 1$ of $L_{p,S}(s,\rho)$, Lemma
  \ref{calculus} interprets the left hand side of (\ref{nice}) as
the leading term at $ s= 0$ of the function $f_\L(\rho\chi_{\rm
cyc}^s)$ with $\L := [R\Gamma_c(U,\T),\zeta_{\Lambda(G)}(\T)] \in
K_1(\La (G),\Sigma_{{\rm ss}-\rho})$. }
\end{remark}

\begin{proof} We note first that if Conjecture \ref{pasc} is valid, then
Remark \ref{colmez} implies that Leopoldt's conjecture is valid
for $E$ and so Lemma \ref{exp-des}(ii) implies that
$R\Gamma_c(U,\T) \in \Sigma_{{\rm ss}-\rho}$ and
$r_G(R\Gamma_c(U,\T))(\rho) = \langle \rho, 1\rangle$ for each
$\rho \in R_{\bc_p}(\bar{G})$.

We now fix $\rho \in R_{\bc_p}(\bar{G})$  and a number field $K$ over which the character
$\rho $ can be realised. We fix an embedding $K \hookrightarrow
 \bc$ and write $\lambda$ for the place of $K$ which is induced by $j$. We set $M := h^0(\Spec E)(1)$ and note that
  $M([\rho ]^*):=M\otimes [\rho ]^*$ is a $K$-motive, where $[\rho]^*$ denotes the dual of the Artin motive corresponding to $\rho.$ 

To evaluate $\zeta_{\Lambda(G)}(\T)^*(\rho)$ we need to make
 Definition \ref{leadingterm-rho} explicit. To do this we use the
observations of \cite[\S1.1, \S1.3]{bufl95} to explicate the
isomorphism $\zeta_{K}(M([\rho ]^*))_{K_\lambda}$ which occurs
in Conjecture \eqref{equivintegrality}. Indeed one has $H^1_f(M) =
\O_E^\times\otimes_\bz\bq$, $H^0_f(M^*(1)) = \bq$, $t_M := E$ and
$H^0_f(M) = H^1_f(M^*(1)) = M_B^+ = 0$ (the latter
 since $E$ is totally real) so that
\[ \bc\otimes_K\Delta_K(M([\rho ]^*)) = \d_\bc((\bq\otimes_\bz\O_E^\times)_{\rho })\d_\bc((\bq)_{\rho })\d_\bc((E)_{\rho })^{-1}\]
and $\zeta_{K}(M([\rho ]^*))_{K_\lambda}$ is equal to the
composite morphism
\begin{align}\label{voi}
 {\bf 1}_{\bc_p} \to &{\bf 1}_{\bc_p}\\
\to
&\d_{\bc_p}((\bq_p\otimes_\bz\O_E^\times)^\rho)\d_{\bc_p}((\bq_p)^\rho)\d_{\bc_p}((\bq_p\otimes_\bq E)^\rho)^{-1}\notag\\
\to
&\d_{\bc_p}(\bc_p\otimes_{K_\lambda}H^2_c(U,M([\rho ]^*)_\lambda))^{-1}\d_{\bc_p}(\bc_p\otimes_{K_\lambda}
H^3_c(U,M([\rho ]^*)_\lambda)\notag\\
\to
&\d_{\bc_p}(\bc_p\otimes_{K_\lambda}R\Gamma_c(U,M([\rho ]^*)_\lambda))^{-1}.\notag
\end{align}
In this displayed formula we have used the following notation: the
first map corresponds to multiplication by
$L_S^*(1,\rho )$; the second map is induced by applying
 $(\bc_p\otimes_{\br, j^{-1}}-)^\rho$ to both the natural isomorphism $E \otimes_\bq \br \cong
\prod_{\Hom(E,\bc)}\br$ and also the exact sequence
\be\label{last-nec} 0 \to \O_E^\times\otimes_\bz \br
\xrightarrow{(\log \circ \sigma)_\sigma} \prod_{\sigma\in
\Hom(E,\bc)}\br \xrightarrow{(x_\sigma)_\sigma \mapsto \sum_\sigma
x_\sigma} \br\to 0;\ee
the third map is induced by Lemma \ref{exp-des}(i) and the inverse
of the isomorphism
\be\label{nec-comp} \prod_{w\in
S_p(E)}U^1_{E_w}\otimes_{\bz_p}\bq_p\xrightarrow{(u_w)_w \mapsto
(\log_p(u_w))_w} \prod_{w \in S_p(E)}E_w \cong \bq_p\otimes_{\bq}E
;\ee
the last map is induced by property h) as described in \S2.1 (with
$R = \bc_p$).

Also, from Lemma \ref{exp-des}(iii) we know that
$\bc_p\otimes_{K_\lambda}t(R\Gamma_c(U,\T)(\rho^*))$ is equal to
the composite
\begin{align}\label{voi2} &\d_{\bc_p}(\bc_p\otimes_{K_\lambda}R\Gamma_c(U,M([\rho ]^*)_\lambda)^{-1} \\
\to
 &\d_{\bc_p}(H^2_c(U,M([\rho ]^*)_\lambda)^{-1}\d_{\bc_p}(H^3_c(U,M([\rho ]^*)_\lambda)\notag\\
 \to &\d_{\bc_p}((\bq_p)^\rho)^{-1}\d_{\bc_p}((\bq_p)^\rho) = {\bf
 1}_{\bc_p}\notag
\end{align}
where the first arrow is induced by property h) and the second by
 Lemma \ref{exp-des}(i) and the map $-c_\gamma^{-1}\log_{p,\gamma,E}$.

Now, after taking account of Lemma \ref{exp-des}(ii),
$\zeta_{\Lambda(G)}(\T)^*(\rho)$ is defined to be $(-1)^{\langle
\rho,1\rangle}$ times the element of $\bc_p^\times$ which
corresponds to the composite of (\ref{voi}) and (\ref{voi2}) (cf.\
Definition \ref{leadingterm-rho}). Thus, after noting that there
is a commutative diagram of the form
\[ \begin{CD} \prod_{w \in S_p(E)}U_{E_w}^1\otimes_{\bz_p}\bq_p @> >>
\cok(\lambda_p)\otimes_{\bz_p}\bq_p\\
@V (\ref{nec-comp}) VV @VV (-1)\times \log_{p,\gamma,E}V\\
 E\otimes_{\bq}\bq_p @> {\rm Tr}_{E/\bq} >> \bq_p\end{CD}\]
where the upper horizontal arrow is the tautological projection,
 the observations made above imply that
\be\label{key} c_\gamma^{\langle \rho, 1\rangle}\cdot
\zeta_{\Lambda(G)}(\T)^*(\rho) =
L_S^*(1,\rho )\cdot\xi\ee
where $\xi$ is the element of
 $\bc_p^\times$ which corresponds to the composite
\begin{align}\label{final-nec} {\bf 1}_{\bc_p} =
&\d_{\bc_p}((\bq_p\otimes_\bz\O_E^\times)^\rho)
\d_{\bc_p}((\bq_p\otimes_\bz\O_E^\times)^\rho)^{-1} \\ \to
&\d_{\bc_p}((\bq_p\otimes_\bq E_0)^\rho)
\d_{\bc_p}((\bq_p\otimes_\bz\O_E^\times)^\rho)^{-1} \notag\\
\to &\d_{\bc_p}((\bq_p\otimes_\bq
E_0)^\rho)\d_{\bc_p}((\bq_p\otimes_\bq E_0)^\rho)^{-1} = {\bf
1}_{\bc_p}\notag
\end{align}
where the first arrow is induced by applying
$\bc_p\otimes_{\br,j^{-1}}-$ to the isomorphism
$\br\otimes_\bz\O_E^\times \cong \br\otimes_\bq E_0$ coming from
the map $(\log \circ \sigma)_\sigma$ in (\ref{last-nec}) and the
second by the isomorphism $\bq_p\otimes_\bz\O_E^\times \cong
\bq_p\otimes_\bq E_0$ coming from the second and third arrows in
(\ref{inf-comp}). (Note that the factor $(-1)^{\langle
\rho,1\rangle}$ in the definition of
$\zeta_{\Lambda(G)}(\T)^*(\rho)$ cancels against the factor $-1$
which occurs in the commutative diagram above, and hence does not
occur in the formula (\ref{key})).

But, upon comparing the definitions of $\mu_\infty$ and $\mu_p$ in
\S\ref{pasc1} with the maps involved in (\ref{final-nec}), one
finds that $\xi$ is equal to
\[ {\rm det}_{\bc_p}((\bc_p\otimes_{\bq_p}\mu_p)\circ (\bc_p\otimes_{\bc,j^{-1}}\mu_\infty)^{-1})^\rho =
 \frac{{\rm det}_{\bc_p}((\bc_p\otimes_{\bq_p}\mu_p)\circ (\bc_p
\otimes_\bq g))^\rho}{({\rm det}_\bc (\mu_\infty\circ (\bc
\otimes_\bq g))^{\rho })}\]
and hence (\ref{key}) implies that
\[ \frac{c_\gamma^{\langle \rho,
1\rangle}\cdot\zeta_{\Lambda(G)}(\T)^*(\rho) }{{\rm det}_{\bc_p}((\bc_p\otimes_{\bq_p}\mu_p)\circ (\bc_p
\otimes_\bq g))^\rho } = \frac{L_S^*(1,\rho )}{{\rm
det}_\bc (\mu_\infty\circ (\bc \otimes_\bq
g))^{\rho }}.\]
The claimed equality (\ref{nice}) now follows immediately upon
comparing this equality to that of Conjecture \ref{pasc}.
\end{proof}

\begin{cor} If Leopoldt's conjecture is valid for $E$ at $p$, then for every $\bq_p$-rational character $\rho$ of $\bar{G}$
there exists a natural number $n_\rho$ such that
\[ (c_\gamma^{\langle \rho, 1\rangle}\cdot\zeta_{\Lambda(G)}(\T)^*(\rho))^{n_\rho} =
 L_{p,S}^*(1,\rho)^{n_\rho}.\]
Further, if $\rho$ is a permutation character, then one can take
$n_\rho = 1$.
\end{cor}

\begin{proof} If $\rho$ is $\bq_p$-rational, then Artin's Induction Theorem implies the existence of a
natural number $n_\rho$ such that in $R_{\bc_p}(G)$ one has
$n_\rho\cdot \rho = \sum_Hn_H\cdot\ind_H^G1_H$ where $H$ runs over
the set of subgroups of $\bar{G}$ and each $n_H$ is an integer (cf.
\cite[chap.\ II, thm.\ 1.2]{tate}). Further, $\rho$ is said to be a
permutation character if and only if there exists such a formula
with $n_\rho = 1$. The stated result thus follows by combining
Theorem \ref{st-thm} with Remark \ref{colmez} and the fact that
each side of (\ref{nice}) is both additive and inductive in
$\rho$.
\end{proof}

\section{The interpolation formula for critical motives}

As a second application of the formalism introduced in \S3, in this section we
prove an interpolation formula for the leading terms (in the sense of
Definition \ref{leadingterm-rho})  of the $p$-adic $L$-functions that Fukaya and Kato conjecture
to exist for any critical motive which has good ordinary reduction at all
places above $p.$ (We recall that a motive $M$ is said to be `critical' if the
map \eqref{alphaM} is bijective). To study these $p$-adic $L$-functions we must  combine
Conjecture \ref{equivintegrality} together with a local analogue of this conjecture (which is also
due to Fukaya and Kato, and is recalled as Conjecture \ref{eps-p} below) and aspects of
Nekov\'a\v{r}'s theory of Selmer complexes and of the theory of $p$-adic height
pairings.

\subsection{Local epsilon isomorphisms}

We fix once and for all 
%the complex period $2\pi
%i,$ i.e.\ a square root of $-1,$ and, for every $l,$ the $l$-adic
the $p$-adic period $t="2\pi i",$ i.e.\ a generator of $\z_p(1).$
Further we write $\epsilon_p(V):=\epsilon(D_{pst}(V))$ for Deligne's $\epsilon$-factor at $p$
where $D_{pst}(V)$ is endowed with the linearized action of the
Weil-group and thereby considered as a representation of the
Weil-Deligne group, see \cite[{\S}3.2]{fukaya-kato}  or
\cite[appendix C]{perrin2000} (Here we suppress the
dependence of the choice of a Haar measure and of $t=2\pi i$ in
the notation. The above choice of $t=(t_n)\in\z_p(1)$ determines a
homomorphism $\psi_p:\qp\to\bar{\qp}^\times$ with
$\ker(\psi_p)=\zp$ sending $\frac{1}{p^n}$ to
$t_n\in\mu_{p^n}).$ Note that $  (B_{dR})^{I_p}=\qnr,$ the completion of the maximal
unramified extension $\mathbb{Q}_p^{nr}$ of $\qp.$  
Let $L$ be any finite extension of $\qp$ and  $V$  any finite dimensional $L$-vector space 
with continuous $G_\qp$-action. We set $\widetilde{L}:= \qnr\otimes_\qp L,$ 
\[ \Gamma^*(-j) := \begin{cases} \Gamma(j)=(j-1)!, &\text{if $j>0$,}\\
 \lim_{s\to j} (s-j)\Gamma(s)=(-1)^{j}((-j)!)^{-1}, &\text{if $j\le 0,$}\end{cases}\]
and $\Gamma_L(V):=\prod_{j \in \bz}\Gamma^*(j)^{-h(-j)},$ where $h(j):=\dim_L gr^j D_{dR}(V).$

In \cite[prop.\
3.3.5]{fukaya-kato}    the existence of a map  
\[\epsilon_{dR}(V):\u_{\widetilde{L}}\to\d_{\widetilde{L}}(V)\d_{\widetilde{L}}(D_{dR}(V))^{-1}\]
with  certain properties was shown. Using this we define an isomorphism
\[\epsilon_{p,L}(V):\u_{\widetilde{L}}\to
\big(\d_L(\r(\qp,V))\d_L(V)\big)_{\widetilde{L}}\]
 as product of $\Gamma_L(V),$ $\eta_l(V)\cdot \overline{(\eta_l(V^*(1))^*)}$  and $\epsilon_{dR}(V).$ 

Now let $T$ be a Galois stable $\O:=\O_L$-lattice of $V$ and set $\widetilde{\O}:=W(\overline{\mathbb{F}_p})\otimes_\zp \O,$ where $W(\overline{\mathbb{F}_p})$ denotes the Wittring of $\overline{\mathbb{F}_p}.$ 

%\begin{conj}[Absolute $\epsilon$-isomorphism]\label{abs-eps-p}
%There exists a (unique) isomorphism
%\[\epsilon_{p,\O}(T): \u_{\widetilde{\O}}\to \big(\d_\O(\r(\qp,T))\d_\O(T)\big)_{\widetilde{\O}}\] with induces $\epsilon_{p,L}(V)$ by base change $L\otimes_\O -.$
%\end{conj}

%This conjecture, which is equivalent to conjecture $C_{EP}(V)$ in
%\cite[III 4.5.4]{fp}, or more precisely its equivariant version
%below  is closely related to the conjecture $\delta_\zp(V)$
%\cite{perrin2000}   via the explicit reciprocity law
%$R\acute{e}c(V),$ which was conjectured by Perrin-Riou and proven
%independently by Benois \cite{benois}, Colmez \cite{colmez1998},
%and Kurihara/Kato/Saito \cite{KKT}. In particular,  the above
%conjecture is known for ordinary cristalline $p$-adic
%representations \cite[1.28,C.2.10]{perrin2000} and for certain
%semi-stable representations, see \cite{berger-tamagawa}.

We define
\[\widetilde{\La}:=W(\overline{\mathbb{F}_p})\kl G\kr=\projlim n \big(
W(\overline{\mathbb{F}_p})\otimes_\zp \zp[G/G_n]\big),\] we assume $L=\qp$ and set as before
$\T:=\La\otimes_\zp T.$  
We write $T(\rho^*)$ for the $\O$-lattice $\rho^*\otimes
T$ of $\rho^*\otimes V,$ which we assume to be de Rham. The following conjecture is a local integrality statement

\begin{conj}[{Fukaya/Kato \cite[conj.\ 3.4.3]{fukaya-kato}, see also
\cite[conj.\ 5.9]{ven-BSD}}]\label{eps-p}
There exists a (unique)  isomorphism in $\C_{\tilde \Lambda}$
\[\epsilon_{p,\La}(\T): \u_{\widetilde{\La}}\to \big(\d_\La(\r(\qp,\T))\d_\La(\T)\big)_{\widetilde{\La}}\] such that for all $\rho:G\to GL_n(\O)\subseteq GL_n(L),$ $L$ a finite extension of $\qp$ with ring of integers $\O,$   we have
\[L^n\otimes_\La \epsilon_{p,\La}(\T)= \epsilon_{p,L}(T(\rho^*)).\]
 \end{conj}

%If $T=T_p\subseteq M_p$ is fixed we also write
%$\epsilon_{p,\La}(M)$ for $\epsilon_{p,\La}(\T).$

\subsection{$p$-adic height pairings}\label{heightpairing}

To prepare for our derivation of the interpolation formula in this
section we discuss certain preliminaries regarding height
pairings.

We fix a continuous finite-dimensional $L$-linear representation
$W$ of $G_\Q$ which satisfies the following `condition of
Dabrowski-Panchiskin':

\begin{itemize}
\item[(DP)] $W$ is de Rham   and there exist a
$G_{\qp}$-subrepresentation $\hat{W }$ of $W$ (restricted to
$G_{\qp}$) such that
$\,D^0_{dR}(\hat{W})=t_p(W):=D_{dR}(W)/D^0_{dR}(W)$.\end{itemize}

Thus we have an exact sequence of $G_{\qp}$-representations
\[ 0 \to \hat{W} \to W \to \tilde{W} \to 0\]
such that $D^0_{dR}(\hat{W})=t_p(\tilde{W})=0$ (cf. \cite[prop.\
1.28]{nek-ht}). Setting $Z:=W^*(1),$ $\hat{Z}:=\tilde{W}^*(1),$
and $ \tilde{Z}:=\hat{W}^*(1),$ we obtain by Kummer duality the
analogous exact sequence
\[ 0 \to \hat{Z} \to Z \to \tilde{Z}\to 0\]
and $Z$ satisfies the condition (DP).

Let $S$ be a finite set of places of $\Q$ containing the places
$S_\infty:=\{\infty\}$ as well as $S_p:=\{p\}$ and such that $W$
(and thus $Z$) are also representations of $G_S(\Q),$ the Galois
group of the maximal outside $S$ unramified extension of $\Q.$
Furthermore, we set $S_f:=S\setminus S_\infty.$

%\subsubsection{Selmer complexes} For $l\neq p$ we define
%$\r_f(\ql,\T)$ as in \eqref{rf-l} and $\r_{/f}(\ql,\T)$ as in
%\eqref{r-/f} with $V$ replaced by $\T,$ see also \eqref{r_flhat}.
%We do {\em not} define $\r_f(\qp,\T)$ since there is in general no
%integral version of $D_{cris}(V).$

The Selmer complex $SC_U(\hat{W},W)$ is by definition the
 mapping fibre \be \label{selmer_U}\xymatrix@C=0.5cm{
    { SC_U(\hat{W},W)\ar[r]^{ } }& {\r(U,W)} \ar[r]^{ } & {\r(\qp,W/\hat{W})\oplus\bigoplus_{S_f\setminus S_p } \r(\ql,W)} \ar[r] &  }\ee
while $SC(\hat{W},W)$ is the mapping fibre
\be \label{selmer}\xymatrix@C=0.5cm{
     {SC(\hat{W},W)}\ar[r]^{ } & {\r(U,W)} \ar[r]^{ } & {\r(\qp,W/\hat{W})\oplus\bigoplus_{S_f\setminus S_p} \r_{/f}(\ql,W)} \ar[r] &  }.\ee

Here $\r_{/f}(\ql,W)$ is defined as mapping cone
\be
\label{r-/f}\xymatrix{
  {\r_f(\ql,W)} \ar[r]^{ } & {\r(\ql, W)} \ar[r]^{ } & {\r_{/f}(\ql,W)} \ar[r]^{ } &     }\ee

%Thus by the octahedral axiom one obtains the distinguished triangles \eqref{sc_u}, \eqref{sc_u-sc} and, using Artin-Verdier/Poitou-Tate duality,
%\be \xymatrix@C=0.5cm{  {\r(U,\T^\vee(1))^\vee} \ar[r]^{} & {SC(\hat{\T},\T)} \ar[r]^{ } & {\r(\qp,\hat{\T})\oplus\bigoplus_{S\setminus (S_p\cup S_\infty)}\r_f(\ql,\T)} \ar[r] &   }. \ee

For any
$G_{\qp}$-representation $V$ and any prime number $\ell$ we define
an element of $L[u]$ by setting
\[ P_\ell(V,u) := P_{L,\ell}(V,u) := \begin{cases} \det_L(1-\varphi_\ell u|
V^{I_\ell}), &\text{if $\ell\neq p$}\\
\det_L(1-\varphi_p u| D_{cris}(V)), &\text{if $\ell = p$,}
\end{cases}\]
where $\varphi_l$ denotes the geometric Frobenius automorphism at
$l.$

The following three conditions are easily seen to be equivalent

\begin{itemize}
\item[(A$_1$)]  $P_l(W,1)P_l(Z,1)\neq0$ for all $l\in S_f\setminus
S_p,$
\item[(A$_2$)] $\H^0(\ql,W)=\H^0(\ql,Z)=0$ for all $l\in
S_f\setminus S_p,$
\item[(A$_3$)] $\r_f(\ql,W)$ is quasi-null for all $l\in
S_f\setminus S_p.$
\end{itemize}

We also consider the following group of conditions

\begin{itemize}
\item[(B$_1$)] $P_p(W,1)P_p(Z,1)\neq0$
\item[(B$_2$)]
$D_{cris}(W)^{\varphi_p-1}=D_{cris}(Z)^{\varphi_p-1}=0$
\item[(B$_3$)] $\H^0(\qp,W)=\H^0(\qp,Z)=0$
\end{itemize}

Then ($B_1$) is equivalent to ($B_2$) and ($B_3$) implies ($B_2$)
by \cite[thm.\ 1.15]{nek-ht}.

Finally we consider the following conditions

\begin{itemize}
\item[(C$_1$)] $P_p(\tilde{W},1)P_p(\tilde{Z},1)\neq0$
\item[(C$_2$)]
$D_{cris}(\tilde{W})^{\varphi_p-1}=D_{cris}(\tilde{Z})^{\varphi_p-1}=0$
\item[(C$_3$)] $\H^0(\qp,\tilde{W})=\H^0(\qp,\tilde{Z})=0.$
\end{itemize}

These conditions are again easily seen to be equivalent (to see
that ($C_2$) is equivalent to ($C_3$) use (loc.\ cit.) and the
fact that $t_p(\tilde{W})=t_p(\tilde{Z})=0).$

\begin{lem}\label{sc-lemma}
Let $X$ denote either $W$ or $Z.$
\begin{itemize}
\item[(i)] If condition ($A_1$) holds, then for every $\ell\in
S_f\setminus S_p$ the following complexes are quasi-null
\[ \r(\ql,X)\cong\r_f(\ql,X)\cong\r_{/f}(\ql,X)\cong 0.\]

\item[(ii)] If condition ($C_1$) is satisfied, then there are
isomorphisms
\[ \r_{/f}(\qp,X)\cong \r(\qp,\tilde{X})\;\;\;\mbox{
and }\;\; \r_f(\qp,X)\cong\r(\qp,\hat{X})\]
in $D^p(L).$

\item[(iii)] If conditions ($A_1$) and ($C_1$) are both satisfied,
then we have a quasi-isomorphism in $D^p(L)$
\[ SC_U(\hat{W},W)\cong\r_f(\Q,W).\]
\end{itemize}
\end{lem}

\begin{proof}
We assume ($A_1$). Then by local duality and the local Euler
characteristic formula it follows immediately that $\r(\ql,X)$ is quasi-null.
The other statements in (i) are obvious. To prove (ii) we assume
(C$_1$). Then, since every complex of vector spaces is  quasi-isomorphic to its cohomology, considered as complex with
zero differential, we have
$\r(\qp,\hat{X})\cong\r_f(\qp,\hat{X})\cong\r_f(\qp,X)$ by
\cite[lem.\ 4.1.7]{fukaya-kato}. Thus the exact triangles
\[ {\r(\qp,\hat{X})} \to {\r(\qp,{X})} \to {\r(\qp,\tilde{X})} \to \]
and
\[ {  \r_f(\qp,X)} \to {\r(\qp,{X})}\to {\r_{/f}(\qp,{X})} \to \]
are naturally isomorphic in $D^p(L).$ Finally, we note that
claim (iii)  follows immediately from
 claims (i) and (ii) and the definitions of $SC_U(\hat{W},W)$ and $ \r_f(\Q,W).$
\end{proof}

Let $M$ be any motive over $\Q,$ $V=M_p$ its $p$-adic realization,
$\rho$ an Artin representation defined over the number field $K$
and $[\rho]$ its corresponding Artin motive. Then, for any place
$\lambda$  of $K$ above $p,$  the $\lambda$-adic realisation \be
W:=N_\lambda=V\otimes_{\qp}[\rho]^*_\lambda\ee of the motive
$N:=M(\rho^*):=M\otimes [\rho]^*$ is an $L:=K_\lambda$-adic
representation. We assume that $V$ (and hence, since $[\rho]^*$ is
pure of weight zero, also $W$) satisfies the condition (DP). We
fix a Galois stable lattice $T$ of $V$ and set
$T_\rho:=T\otimes_{\zp} \O^n,$ a Galois stable lattice in $W$
(where we assume that without loss of generality
$[\rho]^*_\lambda$ is given as $\rho^*: G_\Q\to GL_n(\O)$).
Similarly we fix a $G_\qp$-stable lattice $\hat{T}$ of $\hat{V}$
and take as $\tilde{T}$ the lattice in $\tilde{V}$ induced from
$T.$ Finally we set $\hat{T}_\rho:=\hat{T}\otimes_\zp \O^n$ and
$\tilde{T}_\rho:= \tilde{T}\otimes_\zp \O^n,$ they are Galois
stable $\O$-lattices of $\hat{W}$ and $\tilde{W},$ respectively.

\begin{example} {\em Let $A$ be an abelian variety that is defined over $\bq$ and set $M :=h^1(A)(1)$. If $A$
 has good ordinary reduction at $p$, then $W := N_\lambda$ satisfies the
conditions (DP), ($A_1$), ($B_1$) and ($C_1$) (the last three for
weight reasons -  Weil conjectures, while  more general, the condition (DP) holds for any motive which has good ordinary reduction at $p$ (see \cite{per-ord})).  However, if, for example, 
$A$  is an elliptic curve with (split) multiplicative reduction at
$p$, then the condition ($B_1$) is not satisfied.}
\end{example}

Now we define a  $G_\qp$-stable $\zp$-lattice 
\[ \hat{T}:=T\cap \hat{V}, \]
of $\hat{V}.$ As before let $\T$
denote the big Galois representation $\La\otimes_\zp T$ and put
$\hat{\T}:=\La\otimes_\zp \hat{T}$ similarly. Then $\hat{\T}$ is
$G_\qp$-stable sub-\La-module of $\T.$ In fact, it is a direct
summand of $\T$ and  we have an isomorphism in $\C_{\tilde{\La}}$ \be
\beta:\d_\La(\T^+)_{\tilde{\Lambda}}\cong \d_\La(\hat{\T})_{\tilde{\Lambda}}.\ee
Now the Selmer complexes $ SC_U(\hat{\T},\T)$ and $SC(\hat{\T},\T)$ are defined analogously as for $W$ above.

We point out that $SC_U(\hat{X},X)$ coincides with the Selmer
complex $ \widetilde{\r_f}(X)$ in \cite[(11.3.1.5)]{nek} for
$X\in\{W,Z\}.$ More generally, if we define
\[ \T_{cyc,\rho}:=\La(\Gamma)\otimes T_\rho\]
and similarly $\hat{\T}_{cyc,\rho}$ and $\tilde{\T}_{cyc,\rho},$
then $SC_U(\hat{\T}_{cyc,\rho},\T_{cyc,\rho})$ identifies with the
Selmer complex $ \widetilde{\r_{f, IW}}(\Q_{{cyc}}/\Q,T_\rho)$ defined in
\cite[(8.8.5)]{nek} (with Nekov\'a\v{r}'s local conditions induced by
$T_\ell^+=\hat{\T}_{cyc}(\rho),$ if $\ell=p,$ and $0$ otherwise,
and taking $S_f$ for {\em his} $\Sigma$). Here $\Gamma$ is the
Galois group of the cyclotomic $\zp$-extension $\Q_{cyc}$ of
$\Q.$ Thus we obtain a pairing
\[ h_p(W): \H^1_f(\Q,W)\times\H^1_f(\Q,Z)\to L \]
from \cite[\S 11]{nek} where $h_p(W)$ is called
$\tilde{h}_{\pi,1,1}.$ By \cite[thm.\ 11.3.9]{nek} the pairing
$h_p(W)$ coincides up to sign with the height-pairing constructed
by Schneider \cite{schneider82} (in the case of abelian
varieties), Perrin-Riou \cite{perrin-height} (for semi-stable
representations) and those constructed earlier by Nekov\'a\v{r}
\cite{nek-ht}: see also \S8.1 in (loc.\ cite.), Mazur-Tate
\cite{mazur-tate} and Zarhin \cite{zarhin} for alternative definitions of related
height pairings.

It follows from the construction of Nekov\'a\v{r}'s height pairing (cf. \cite[the
sentence after
 (11.1.3.2)]{nek}) that the induced map
\be\label{ad} \mathrm{ad}(h_p(W)):\H^1_f(\Q,W)\to\H^1_f(\Q,Z)^*\ee
is equal to the composite
\begin{multline}\label{ad-bock} {\H^1_f(\Q,W)}\cong\H^1(SC_U(\hat{W},W))
 \xrightarrow{\B} {\H^2 (SC_U(\hat{W},W))}\\\cong\H^2_f(\Q,W)\cong
\H^1_f(\Q,Z)^*\end{multline}
where the first and third maps are by Lemma \ref{sc-lemma}(iii),
$\B$ denotes the Bockstein morphism for
 $SC_U(\hat{\T}_{cyc,\rho},\T_{cyc,\rho})$ and the last map is by global
duality.

 For the evaluation at representations in the next subsection we  need  the following descent properties.  Let $\Upsilon$ be the set of all $l\neq p$ such that the ramification index of $l$ in $F_\infty/\Q$ is infinite. Note that $\Upsilon$ is empty if $G$ has a commutative open subgroup.

 \begin{prop}\label{descent-prop}\cite[prop.\ 1.6.5]{fukaya-kato}
 With the above notation  we have canonical isomorphisms (for all $l$)
 \begin{align*}
L^n\otimes_{\La,\rho} \r(U,\T)&\cong \r(U,W), &  L^n\otimes_{\La,\rho} \r_{c}(U,\T)&\cong \r_{c}(U,W), \\
  L^n\otimes_{\La,\rho} \r(\ql,\T)&\cong \r(\ql,W), &
L^n\otimes_{\La,\rho}SC_U(\hat{\T},\T)&\cong SC_U(\hat{W},W).
 \end{align*}
 For $l\not\in\Upsilon\cup S_p$ we also have: $L^n\otimes_{\La,\rho} \r_f(\ql,\T)\cong
 \r_f(\ql,W).$
  \end{prop}

But note that the complex $\r_f(\ql,\hat{\T})$ for $l\in\Upsilon$
and thus $SC(\hat{\T},\T)$ does {\em not} descend like this in
general. Instead, according to \cite[prop.\ 4.2.17]{fukaya-kato}
one has a distinguished triangle \be
 \xymatrix@C=0.5cm{
  L^n\otimes_{\La,\rho}SC(\hat{\T},\T)\ar[rr]^{ } && SC(\hat{W},W) \ar[rr]^{ } && {\bigoplus_{l\in \Upsilon}\r_f(\ql,W)} \ar[r] &   }.\ee

\subsection{The interpolation formula}\label{tif}

In this section we assume that the motive $N :=M(\rho^*)$ is
critical. Then, assuming \cite[conj.\ 3.3]{ven-BSD} of Fontaine
and Perrin-Riou it follows that the motivic cohomology groups
\begin{itemize}
\item[($D_1$)] $\H_f^0(N)=\H^0_f(N^*(1))=0$
\end{itemize}
both vanish. Assuming also the conjecture \cite[conj.\
3.6]{ven-BSD} on $p$-adic regulator maps, this is equivalent to
the condition
\begin{itemize}
\item[($D_2$)] $\H^0_f(\Q,W)=\H^0_f(\Q,Z)=0.$
\end{itemize}

%\begin{example} {\em If $A$ is an abelian variety over $\bq$, then the conditions ($D_1$) and
%($D_2$) are both satisfied for the motive
%$M=h^1(A)(1)$.}\end{example}

We also consider the condition
\begin{itemize}
\item[(F)] The pairing $h_p(W)$ is non-degenerate.
\end{itemize}

\begin{example}
{\em If $A$ is an abelian variety over $\bq$, then the conditions ($D_1$) and
($D_2$) are both satisfied for the motive
$M=h^1(A)(1)$.
However, very little is known about the non-degeneracy of the p-adic height pairing in the ordinary case. Indeed, as far as we are aware, the only theoretical evidence
for non-degeneracy is an old result of Bertrand  \cite{ber} that for an elliptic curve with complex multiplication, the height of a point of 
infinite order is non-zero (even this is unknown in the non CM case).
Computationally, there has been a lot of work done recently by Stein and 
  Wuthrich \cite{wuth}. We are grateful  to J. Coates,  P. Schneider and C. Wuthrich for providing us with these examples.}
\end{example}

\begin{prop} Assume that the conditions ($A_1$), ($C_1$) and ($D_2$) are
satisfied. Then 
\begin{itemize}
\item[(i)] $SC_U(\hat{\T}_{cyc,\rho},\T_{cyc,\rho})\in \Sigma_{ss}$   if and only if the condition
(F) holds.
\item[(i)] if (F) holds, we have $r_\Gamma(SC_U(\hat{\T}_{cyc,\rho},\T_{cyc,\rho}))=\dim_L \H^1_f(\bq,W).$
\end{itemize}
\end{prop}

\begin{proof}
By assumption ($D_2$) implies that $\H^i(SC_U(\hat{W},W))=0$ for
$i\neq 1,2.$ Thus the claim follows from the fact that \eqref{ad}
and \eqref{ad-bock} coincide.
\end{proof}

Now let $F_\infty$ be as before a $p$-adic Lie extension of $\Q$ which contains
$\Q_{cyc}$ and   $G$   its Galois group (with  quotient $\Gamma$). We
set $\La := \La(G).$ Since by \cite[4.1.4(2)]{fukaya-kato} we have
a canonical identification \be
(\La_\O(\Gamma)\otimes_\O\O^n)\otimes_{\La(G)}^\mathbb{L}SC_U(\hat{\T},\T)\cong
SC_U(\hat{\T}_{cyc,\rho},\T_{cyc,\rho})\ee we conclude that
$SC_U(\hat{\T},\T)\in\Sigma_{{\rm ss}-\rho}$ if and only (F) holds
for $W.$

In \cite{fukaya-kato} two $p$-adic $L$-functions \be
\L_U:=\L_U(M):\u_\La\to \d_\La(SC_U(\hat{\T},\T))\ee and \be
\L:=\L(M):\u_\La\to \d_\La(SC(\hat{\T},\T))\ee were defined, modulo the validity of
 Conjecture \ref{equivintegrality} and Conjecture \ref{eps-p}.   Recall the definition of $\Sigma_{SC_U}$ or $\Sigma_{SC}$ at the end of section \ref{prelim} where we abbreviate $SC_U=SC_U(\hat{\T},\T)$ and $SC=SC(\hat{\T},\T).$  Now $\L_U$ and $\L$ induce classes $[SC_U,\L_U]$ and $[SC,\L]$ in
$K_1(\La(G),\Sigma_{SC_U})$ and $K_1(\La(G),\Sigma_{SC}),$ which for simplicity we again call $\L_U$ and $\L,$ 
respectively. However, we note that for comparison purposes, it would perhaps
be more convenient to define both $\L$ and $\L_U$ as elements of the group $K_1(\La(G),\Sigma'_{SC})$ where $\Sigma'_{SC}$ is the smallest full subcategory of $C^p(\La)$
which satisfies  conditions (i)-(iii) and (iv') from section \ref{prelim} and
contains $SC_U$ as well as the Euler factors $\r_f(\ql,\hat{\T})$ for all $l\in S_f\setminus S_p.$

\begin{thm}\label{int-thm} We assume that $W$ satisfies the conditions $(A_1)$, $(B_1)$, $(C_1)$, $(D_1)$ and $(F)$ and that the
 isomorphisms $\zeta_{\La}(M)$ and $\epsilon_{p,\La}(\hat{\T})$ that are described in Conjectures \ref{equivintegrality} and \ref{eps-p}  both exist. Then both  $SC_U(\hat{\T},\T)$ and $ SC(\hat{\T},\T)$ belong to  $ \Sigma_{\rm ss},$ $r:=r_G(SC_U(\hat{\T},\T))(\rho)=r_G( SC(\hat{\T},\T))( \rho)= \dim_L \H^1_f(\bq,W)$ and   the
leading term $\L^*(\rho)$ (respectively $\L^*_U(\rho)$)
is equal to \be\label{interpolation} (-1)^r
\frac{L_{K,B}^*(M(\rho^*))}{\Omega_\infty(M(\rho^*))R_\infty(M(\rho^*))}\cdot\Omega_p(M(\rho^*))R_p(M(\rho^*))\cdot\Gamma(\hat{V})^{-1}
\cdot \frac{P_{L,p}(\hat{W}^*(1),1)}{P_{L,p}(\hat{W},1) },\ee
where $L_{K,B}^*(M(\rho^*))$ is the leading coefficient of
the    complex $L$-function of $N,$ truncated by removing Euler factors
for all primes in  $B:=\Upsilon\cup S_p$ (respectively $B:=S\setminus S_\infty$), and
the  regulators $ R_\infty(M(\rho^*))$ and $ R_p(M(\rho^*))$ and periods
$\Omega_\infty(M(\rho^*))$ and $ \Omega_p(M(\rho^*))$ are as defined in the course
of the proof given below.
\end{thm}

\begin{remark} {\em The formulas of Theorem \ref{int-thm} generalize the formulas obtained by Perrin-Riou in  
\cite[4.2.2,4.3.6]{perrin2000}. Further, by slightly altering the definition of the  complex $L$-function  an analogous formula can be proved even in the case that ($B_1$) is not satisfied. Indeed, if ($B_1$)  fails, then one can have $P_{L,p}({W},0)=0$ and so the order of vanishing of $L_{K,B}(M(\rho^*),s)$ and $L_{K}(M(\rho^*),s)$  may differ. However, to avoid this problem, in formula \eqref{interpolation}  one need only replace $P_{L,p}(\hat{W},1)$ by the leading coefficient of $P_{L,p}(\hat{W},p^s)$ at $s=0,$ or equivalently one can
replace the term  $\frac{L_{K,B}^*(M(\rho^*))}{P_{L,p}(\hat{W},1)}$ by  $\frac{L_{K,B\setminus\{p\}}^*(M(\rho^*))}{\{P_{L,p}({W},u)^{-1}P_{L,p}(\hat{W},u)\}_{u=1}}.$ }  
\end{remark}

\begin{proof}
 This proof is closely modelled on that of \cite[thm.\ 4.2.26]{fukaya-kato} (as
amplified in \cite[proof of thm.\
6.4]{ven-BSD}). At the outset we let $\gamma=(\gamma_i)_i$ and $\delta=(\delta_i)_i$ be a choice of `good bases (in the sense of   \cite[4.2.24]{fukaya-kato}) of $M_B^+$ and $t_M$ for
$\hat{\T}$ and let $\gamma'$ and $\delta'$ be the induced
$K$-bases of $N_B^+$ and $t_N,$
 respectively. This induces a map \be can_{\gamma'.\delta'}:\u_K\to
 \d_K(N_B^+)\d_K(t_N)^{-1}.\ee
Furthermore, let $P^\vee=(P_1^\vee,\ldots,P_{d(N)}^\vee)$ and
$P=(P_1,\ldots,P_{d(N)})$ be $K$-bases of $\H^1_f(N)$ and
$\H^1_f(N^*(1)),$ respectively. Setting $P^d=(P_1^d,\ldots,
P^d_{d(N)})$ the dual basis of $P$ we obtain similarly \be
can_{P^\vee,P}:\u_K\to\d_K(\H^1_f(N))\d_K(\H^1_f(N^*(1))^*)^{-1}.\ee
Then $can:=can_{\gamma'.\delta'}\cdot can_{P^\vee,P}$ is an
isomorphism \be can:\u_K\to\Delta_K(N)=
 \d_K(N_B^+)\d_K(t_N)^{-1}\d_K(\H^1_f(N))\d_K(\H^1_f(N^*(1))^*)^{-1}.\ee
We fix an embedding of $K$ into $\bbC.$
 Now let $\Omega_\infty(N)$ and $R_\infty(N)$ denote the
 determinant of the isomorphism
 \be\alpha_N: (N_B^+)_{\bbC}\to (t_N)_{\bbC}\ee
 with respect to $\gamma'$ and $\delta',$ and the inverse of
 \be h_\infty:\big(\H^1_f(N^*(1))^*\big)_{\bbC} \to \H^1_f(N)_{\bbC}\ee
 with respect to $P^d$ and $P,$ respectively. In other words, we
 have
 \[ \Omega_\infty(N): \,\,\,{\u_\bbC} \xrightarrow{(can_{\gamma',\delta'})_\bbC}
 {\d_K(N_B^+)_\bbC\d_K(t_N)_\bbC^{-1} } \xrightarrow{\d(\alpha_N)\cdot \id} {\u_\bbC}
 \]
and
\[ R_\infty(N): \,\,\,{\u_\bbC} \xrightarrow{(can_{P^\vee,P})_\bbC}
{\d_K(\H^1_f(N))_\bbC\d_K(\H^1_f(N^*(1))^*)^{-1}_\bbC}
\xrightarrow{\id\cdot\d(h_\infty)^{-1}} {\u_\bbC}.\]
Note that we obtain
\[ \Omega_\infty(N)R_\infty(N): \,\,\, {\u_\bbC} \xrightarrow{can} {\Delta_K(N)_\bbC }\xrightarrow{(\vartheta_\infty)_\bbC}
  {\u_\bbC}.\]
Upon comparing this with
\[ L^*_K(M):\,\,\,\, {\u_\bbC} \xrightarrow{\zeta_K(N)_\bbC} {\Delta_K(N)_\bbC }\xrightarrow{(\vartheta_\infty)_\bbC}
  \u_\bbC.   \]
we deduce that $\zeta_K(N):\u_K\to \Delta_K(N)$ coincides with
\[ \frac{L^*_K(M)}{\Omega_\infty(N)R_\infty(N)}\cdot can:\u_K\to
   \Delta_K(N).\]

Also, $\L_U^*(\rho)$ is
defined (in Definition \ref{leadingterm-rho}) to be the isomorphism
\begin{multline*} {\u_{\tilde{L}}} \xrightarrow{\zeta_\La(M)(\rho)_{\tilde{L}}}
{\d_L(\r_c(U,W))^{-1}_{\tilde{L}}}
\xrightarrow{\beta(\rho)\epsilon(\hat{\T})^{-1}(\rho)}\\
{\d_L(SC_U(\hat{W},W))^{-1}_{\tilde{L}}} \xrightarrow{t(SC_U(\rho^*))_{\tilde{L}}}
\u_{\tilde{L}}\end{multline*} where we set $\zeta_\La(M)(\rho):={L}^n\otimes_\La \zeta_\La(N),$ $\beta(\rho):=L^n\otimes_\La \beta$  and $\epsilon(\hat{\T})(\rho):=L^n\otimes_\La \epsilon_{p,\La}(\hat{\T}).$
Using that $\zeta_\La(M)(\rho) $ equals
\[ {\u_{\tilde{L}}} \xrightarrow{\zeta_K(N)_{\tilde{L}}} {\Delta_K(N)_{\tilde{L}}} \xrightarrow{\vartheta_\lambda}
{\d_L(\r_c(U,W))^{-1}_{\tilde{L}}}\]
by \cite[conj.\ 2.3.2]{fukaya-kato} (or \cite[conj.\ 3.7, Conj.
4.1]{ven-BSD}), one sees easily that $\L_U^*(\rho)$ is the product
of the following terms \eqref{1}-\eqref{7}:

%\begin{itemize}
  \be\label{1}\frac{L_K^*(N)}{ \Omega_\infty(M(\rho^*))
R_\infty(N)},\ee
%\item[(2)-(7)] as in (loc.\ cit.) before lemma 6.3,
 \be\label{2}\Gamma_L(\hat{W})^{-1}=\Gamma_\qp(\hat{V})^{-1},\ee
 \begin{multline}\label{3}
\Omega_p(M(\rho^*)):\xymatrix{
   {\d_L(\hat{W})_{\tilde{L}}}\ar[r]^{\cdot\;\epsilon_{dR}(\hat{W})^{-1}} &{\d_L(D_{dR}(\hat{W}))_{\tilde{L}}} \ar[r]^{\d(g_{dR}^t)}&{\d_K(t_{M(\rho^*)})_{\tilde{L}}}\ar[r]^<(0.2){\cdot\; can_{\gamma,\delta}}&}\\\xymatrix{
   {\d_K\big(M(\rho^*)_B^+\big)_{\tilde{L}}} \ar[r]^<(0.3){\d(g_\lambda^+)}   & {\d_L(W^+)_{\tilde{L}}}\ar[r]^{\beta(\rho)}&  { \d_L(\hat{W})_{\tilde{L}}} } 
 \end{multline} %which is, by definition, the composite
 
where we apply Remark \ref{inverse} to obtain an automorphism of
$\u_{\tilde{L}}.$
%\footnote{Using Remark \ref{inverse}(i) it is easy to see that this amounts to taking the
%product of the following isomorphisms and identifying the target with $\u$ afterwards
%\begin{align*}
%&\xymatrix{{\u}\ar[rr]^<(0.2){can_{\gamma,\delta}} &   &
%{\d\big((M(\rho^*)_B^+)_L\big)\d(t_{M(\rho^*)})^{-1},} } &
%&\xymatrix{ {\u}\ar[rr]^<(0.25){\id_{-}\cdot{\d(g_\lambda^+)}  } & &{ \d(W^+) \d\big((M(\rho^*)_B^+)_L\big)^{-1},}}\\
%&\xymatrix{ {\u}\ar[rr]^<(0.25){\id_{-}\cdot{\d(g_{dR}^t)}} &   &
%{\d(D_{dR}(\hat{W}))^{-1}} \d(t_{M(\rho^*)}),} &
%&\xymatrix{ {\u}\ar[rr]^<(0.3){\epsilon_{dR}(\hat{W})^{-1}} &   & {{ \d(\hat{W})^{-1}}\d(D_{dR}(\hat{W}))},   }\\
%&\xymatrix{ {\u}\ar[rr]^<(0.3){\id_{-}\cdot\beta(\rho)} &   &{
%\d(\hat{W})}{ \d(W^+)^{-1},} }
%\end{align*}
%where the identity maps are those of $\d\big((M(\rho^*)_B^+)_L\big)^{-1},$ $\d(D_{dR}(\hat{W}))^{-1}$ and $\d(W^+)^{-1},$ respectively.}
  \be\label{4}\prod_{S\setminus\{p,\infty\}} P_{L,l}(W,1):\xymatrix{
  {\u_L} \ar[rr]^<(0.3){\prod\eta_l(W)}   && {\prod \d_L(\r_f(\ql,W))} \ar[r]^<(0.4){acyc} & {\u_L }  }\ee where the first map comes from the trivialization by the identity and the second from the acyclicity,\\
\begin{multline}\label{5}
  \{P_{L,p}(W,u)P_{L,p}(\hat{W},u)^{-1}\}_{u=1}:\\\xymatrix{
  {\u_L} \ar[rr]^<(0.25){\eta_W\cdot\eta_{\hat{W}}^{-1}} &&
  {\d_L(\r_f(\qp,W)) \d_L(\r(\qp,\hat{W}))^{-1}} \ar[rr]^<(0.4){quasi} &&
   {\u_L , }   } \end{multline} where we use that $t_p(W)=D_{dR}(\hat{W})=t_p(\hat{W})$ and the quasi-isomorphism mentioned in Lemma \ref{sc-lemma}(ii), \\
\be \label{6} P_{L,p}(\hat{W}^*(1),1):\xymatrix{
  {\u_L} \ar[rr]^<(0.2){\overline{(\eta_{\hat{W}^*(1)})^*}} && {\d_L(\r_f(\qp,\hat{W}^*(1)))} \ar[r]^<(0.4){acyc} & {\u_L,} }\ee where we use that $t_p(\hat{W}^*(1))=D_{dR}^0(\hat{W})=0,$ and 
 the determinant over $L$ of
    the isomorphism $ad(h_p(W))$ with respect to the bases
    $P^\vee$ and $P$:
\begin{multline}\label{7}
 R_p(N):\,\,\,\, \u_L\xrightarrow{(can_{P^\vee,P})_L} {\d_K(\H^1_f(N))_L\d_K(\H^1_f(N^*(1))^*)_L^{-1}} \xrightarrow{\cong}
 \\ {\d_K(\H^1_f(\Q,W))_L\d_K(\H^1_f(\Q,Z)^*)_L^{-1}} \xrightarrow{h_p(W)} {\u_L}. \end{multline}
 Indeed, as remarked above $\zeta_\La(M)(\rho)$ decomposes up to the comparison isomorphism $\d(g_\lambda^+),$ which contributes to factor \eqref{3}, into $\zeta_K(N)_L$ and $\vartheta_\lambda$. While $\zeta_K(N)_L$ gives the full factor \eqref{1} and contributes with $can_{\gamma,\delta}$ and $can_{P^\vee,P}$ to the factors \eqref{3} and \eqref{7}, respectively, the second part $\vartheta_\lambda$ gives the full factor \eqref{4}, the half factor \eqref{5} in form of $\eta(W)$ and contributes $\d(g_{dR}^+)$ to factor \eqref{3}. Further, $\beta(\rho)$ contributes to factor \eqref{3}, while according to   \cite[\S 3.3]{fukaya-kato} $\epsilon(\hat{\T})^{-1}(\rho)=\epsilon_{p,L}(\hat{W})^{-1}$  gives the full factors \eqref{2} and \eqref{6}, the other half of \eqref{5} in form of $\eta_{\hat{W}}^{-1}$ and adds $\epsilon_{dR}(\hat{W})$ to factor \eqref{3}. Finally,  we had observed at the end of  \S \ref{heightpairing} that $t(SC_U(\rho^*))$ equals $h_p(W).$ 

    In order to derive the corresponding formula for $\L^*(\rho)$
    we use the   exact triangle
\[ SC_U(\hat{W},W) \to L^n\otimes^\mathbb{L}_\La SC(\hat{\T},\T) \to
{\bigoplus_{S_f\setminus(S_p\cup\Upsilon)}\r_f(\ql,W)} \to \]
and  the equality \[ L_{K,\Upsilon'}^*(N)=L_{K,B'}^*(N) \prod_{B\setminus\Upsilon}
P_{L,l}(W,1)^{-1}\]
with $\Upsilon'=\Upsilon\cup \{p\}$ and $B'=S\setminus S_\infty.$
\end{proof}
 
\begin{example}
 {\em Let $E$ be an elliptic curve defined over $\bq,$ $M=h^1(E)(1)$ and $F_\infty:=\bq (E(p))$ where $E(p)\subseteq E(\overline{\bq}]$ denotes the group of $p$-power torsion points of $E.$ It is conjectured that $SC_U(\hat{\T},\T)$ always belongs to $\Sigma_{S^*}$ (cf. \cite[conj.\ 5.1]{cfksv} and \cite[4.3.5 and
prop.\ 4.3.7]{fukaya-kato}). As was shown in \cite{fukaya-kato} the existence of $\L(M)$ implies the existence of $\L_E\in K_1(\La(G)_{S^*})$ as conjectured in \cite[conj.\ 5.7]{cfksv} to exist  with a prediscribed interpolation property $\L_E(\rho)$ for $r_G(SC(\hat{\T},\T))(\rho)=0.$ Now the above formula \eqref{interpolation} predicts   the leading terms $\L_E^*(\rho)$ for all   Artin characters with non-degenerate archimedean and $ p$-adic height pairing. } 
\end{example}

\Addresses
\end{document}